\documentclass[12pt]{article}
\usepackage{amsmath}
\usepackage{amssymb}
\usepackage{amsthm}
\usepackage{amscd}
\usepackage{amsxtra}
\usepackage{verbatim}
\usepackage{xcolor}
\usepackage{color}
\usepackage{enumerate}
\usepackage{mathrsfs}
\usepackage{array,arydshln}
\usepackage{bm}
\allowdisplaybreaks
\usepackage{amsfonts}
\usepackage{mathrsfs}

\numberwithin{equation}{section}
\definecolor{dblue}{rgb}{0,0,0.45}
\definecolor{red}{rgb}{0.7,0,0}

 \RequirePackage{geometry}
\newtheorem{theorem}{Theorem}[section]
\newtheorem{lemma}[theorem]{Lemma}
\newtheorem*{lemma*}{Lemma}

\newtheorem{proposition}[theorem]{Proposition}

\theoremstyle{definition}

\newtheorem{remark}[theorem]{Remark}
\newtheorem{definition}[theorem]{Definition}

\theoremstyle{remark}

\def\para#1{\vskip .4\baselineskip\noindent{\bf #1}}

\date{\today}

\begin{document}
\title{Averaging principle for fast-slow system driven by 
	mixed fractional Brownian rough path \author{ Bin \textsc{Pei}, Yuzuru \textsc{Inahama} and Yong \textsc{Xu}}}
\maketitle

\begin{abstract}
This paper is devoted to studying the averaging principle for fast-slow system of rough differential equations driven by mixed fractional Brownian rough path.  The fast component is driven by Brownian motion, while the slow component is driven by  fractional Brownian motion with Hurst index $H ~(1/3 <  H\leq 1/2)$.  Combining the fractional calculus approach to rough path theory 
and Khasminskii's classical time discretization method,
we prove that the slow component strongly converges to the solution of the corresponding averaged equation in the $L^1$-sense.
The averaging principle for a fast-slow system 
in the framework of rough path theory seems new.
\end{abstract}

\section{Introduction}
Let $(\Omega, \mathcal{F}, \{ \mathcal{F}_t\}_{t\geq 0}, P)$ be a complete probability space, $W=(W_t)_{t\geq 0}$ be a standard $d'$-dimensional Brownian motion (Bm), $B=(B_t)_{t\geq 0}$ be an $d$-dimensional fractional Brownian motion (fBm) with Hurst index $H\in (\frac13, \frac12]$, that is a collection of
centered, independent Gaussian processes, independent of $W$ as well, with covariance
function
\begin{eqnarray*}
	R_{H}(t, s)=\frac{1}{2}\left(t^{2 H}+s^{2 H}-|t-s|^{2 H}\right)I_d, \quad s, t \geq 0,
\end{eqnarray*}
where $I_d$ is the identity matrix of size $d$.
The Kolmogorov theorem entails that $\mathrm{fBm}$ has a modification with $\beta$-H\"older sample paths for any $\beta<H$. For each $t\geq 0$, we denote by  $\mathcal{F}_t$ the $\sigma$-field generated by the random variables $\{B_s,W_s,s \in[0,t]\}$ and all $P$-null sets. The expectation with respect to $P$ is denoted by $\mathbb{E}$. 
In addition to the natural filtration $\{\mathcal{F}_t,t\geq 0\}$, we will consider a larger filtration $\{\mathcal{G}_t,t\geq 0\}$ such that $\{\mathcal{G}_t\}$ is right-continuous and $\{\mathcal{G}_0\}$ contains the  $P$-null sets, so that $B$ are $\mathcal{G}_0$-measurable, and $W$ is a $\{\mathcal{G}_t\}$-Brownian motion.	

In what follows, we will denote by  $C^{k}_b(\mathbb{R}^m;\mathbb{R}^n)$ the set of functions $f:\mathbb{R}^m\rightarrow \mathbb{R}^n$ which
are bounded, $k$-times continuously differentiable 
 with bounded derivatives of order up to $k$
 (in symbols $\nabla f, \nabla^2f,\ldots,\nabla^k f$).
We denote by $C^{k,\lambda}_{b}(\mathbb{R}^m;\mathbb{R}^n)$
the set of $f \in C^{k}_b(\mathbb{R}^m;\mathbb{R}^n)$ whose
 $k$th derivative is
uniformly H\"older continuous of order $\lambda\in(0,1]$. 
The time interval will be $[0,T]$ for arbitrary $T>0$.
The transpose of a vector $A$ will be denoted by $A^{{\rm T}}$.

We firstly consider the differential equation driven by both Bm
and  fBm of the type
\begin{equation}\label{msde}
u_t=u_0+\int_{0}^{t}a(u_s) \mathrm{d}  s+\int_{0}^{t}b(u_s)\mathrm{d}  B_s+\int_{0}^{t}c(u_s) \mathrm{d}  W_s,
\end{equation}
with $u_0 \in \mathbb{R}^e$ which is arbitrary and non-random but fixed.
Since the fBm is neither Markov nor semimartingale if $H\neq \frac12$, 
we cannot use the classical It\^o theory 
to solve (\ref{msde})  unless $b\equiv 0$.  
Lyons \cite{lyons1998} made a breakthrough by inventing rough path theory, which enabled us to do
 pathwise study of stochastic differential equations (SDEs) as above.
We will show that  (\ref{msde}) can be understood as a rough differential 
equation (RDE) and possesses a unique global solution
if the coefficients $a:\mathbb{R}^e \rightarrow \mathbb{R}^e$, $b:\mathbb{R}^e \rightarrow \mathbb{R}^{e}\otimes \mathbb{R}^{d}$ and $c:\mathbb{R}^e \rightarrow \mathbb{R}^{e}\otimes \mathbb{R}^{d'}$ satisfy suitable regularity assumptions. 

The driving rough path of RDE (\ref{msde}) is a natural 
rough path lift of $(B_t, W_t)_{0 \le t \le T}$,
which is formally given by $\mathbf{Z} =(Z, Z^2)$, where
\begin{equation}\label{RPZ}
Z_{t}=(B_{t},W_{t})^{{\rm T}}\quad {\rm and} \quad Z^2_{st}=\bigg(\begin{array}{ll}B^2_{st}& \int_{s}^{t}
(B_{u} -B_s)\otimes \mathrm{d}W_u
\\ \int_{s}^{t} (W_{u} -W_s) \otimes \mathrm{d}B_u
& W^2_{st} 
\end{array}\bigg).
\end{equation}
For every $\beta \in (\frac13, H)$, this lift $\mathbf{Z}$ exists 
and is a $\beta$-H\"older (weakly) geometric rough path almost surely (see Section 3 for details). 
We call it mixed fractional Brownian rough path.
Since the $W$-component of $\mathbf{Z}$ is Stratonovich-type
Brownian rough path, the last term in (\ref{msde}) is 
something like a Stratonovich integral.

In the case $H\in (\frac12,1)$, it is well-known that Young integral, which is essentially a generalized Riemann-Stieltjes integral,
could be a good choice to give meaning to the integral with respect to fBm \cite{Rascanu2002,Zahle1998}. 
When $c \equiv0$ in (\ref{msde}), the theory of SDEs driven only by fBm with $H\in (\frac12,1)$ was initiated by Lyons \cite{lyons1994}
and has been well developed by many authors
especially on the existence and uniqueness of pathwise solutions. 
For example,
using the fractional calculus introduced by Z\"ahle \cite{Zahle1998}, Nualart and R\u{a}\c{s}canu \cite{Rascanu2002} derive very weak conditions for the general case where $a,b$ are functions of $(t,u_t)$, in particular
$b$ need to be only $C^1$ with bounded and H\"older continuous first derivative, to ensure the
existence and uniqueness of the solution in the space of H\"older continuous functions.  
When $c$ does not vanish identically, 
there are only a few results devoted to such mixed equations. The main difficulty when
considering (\ref{msde}) lies in the fact that both stochastic integrals are dealt in
different ways. The integral with respect to the Bm is understood as an It\^o integral, while the integral with respect to the fBm has to be understood in the pathwise Young sense. 
Kubilius \cite{Kubilius2002} studies SDEs driven by both
fBm and Bm, in the one-dimensional
case, with no drift term. Guerra and Nualart \cite{Guerra2008} combine the pathwise approach (generalized Riemann-Stieltjes integral) with the It\^o stochastic calculus to prove an existence and uniqueness theorem for 
multidimensional, time-dependent SDEs driven
simultaneously by a multidimensional fBm with Hurst
parameter $H > 1/2$ and a multidimensional Bm.

Let us get back to the case $H\in (\frac13, \frac12]$.
When $c\equiv0$ in (\ref{msde}), the existence and uniqueness result was generalized by
Lyons' seminal paper \cite{lyons1998}. 
In this paper he established rough path theory to define 
the integral with respect to $\beta$-H\"older rough path 
 ($\frac13 <\beta < \frac12$).
This theory is basically deterministic (see several monographs \cite{FrizV2010, FrizH2017, Lyons2002, Lyons2007}). 
 Coutin and Qian \cite{Coutin2002} proved that fBm admits a 
natural rough path lift in $1/\beta$-variation 
rough path topology, which was later 
improved to $\beta$-H\"older rough path topology ($\frac13 <\beta <H$).
These results enabled us to study   (\ref{msde}) 
via rough paths when  $c\equiv0$.

There are other formulations of rough path theory.
Gubinelli \cite{Gubinelli2004} established an alternative theory of controlled rough paths to generalize the concept of integration and differential equations with H\"older exponent greater than $\frac13$. 
Following Z\"ahle \cite{Zahle1998}, 
Hu and Nualart \cite{HuN2009} developed another approach
to rough path theory by using fractional calculus.

Like the case of usual SDEs, the condition on the drift $a$
can be weaker than one on the diffusion coefficient $b$.
A well-known result states that a unique global solution exists when $a \in C^1_b$ and $b \in C_b^3$
(see \cite{Garrido2018} for instance).
Riedel and Scheutzow \cite{Riedel2017} solved RDEs with
unbounded drift term. In this work, 
$a$ is allowed to grow at most linearly,  while  $b \in C_b^4$
 (see Proposition \ref{uniqpro} below).

When $c$ does not vanish identically,  however,
 much less is known about (\ref{msde}). 
The most relevant result in our case is
Diehl, Oberhauser and Riedel \cite{Diehl2015}, in which 
the authors gave a meaning to differential equations driven by
 a deterministic rough path and  Brownian rough path. 
 In contrast to the RDE in \cite{Diehl2015}, 
 the trajectories of $B$ and $W$ in (\ref{msde}) are both stochastic. 

Now, we summerize the  
 existence and uniqueness result for (\ref{msde}). 
Needless to say, the unique solution does not depend on the choice of $\beta \in (\tfrac13, H)$.
\begin{proposition}\label{uniqpro}
Let $\tfrac13 < \beta < H \le \tfrac12$ and write
$b=(b_1,\ldots,b_d)$ and $c=(c_1,\ldots, c_{d'})$.
Assume either one of the following two conditions on the coefficients  of RDE (\ref{msde}).
\begin{enumerate}
\item 
$a$ is a locally Lipschitz continuous vector field with at most linear growth on $\mathbb{R}^{e}$ and $b_i,c_j \in C_b^{4}(\mathbb{R}^{e}, \mathbb{R}^{e})~(1\leq i \leq d, 1\leq j \leq d')$. 
\item $a\in C_b^{1}(\mathbb{R}^{e}, \mathbb{R}^{e})$ and $b_i,c_j \in C_b^{3}(\mathbb{R}^{e}, \mathbb{R}^{e})~(1\leq i \leq d, 1\leq j \leq d')$. 
\end{enumerate}
Then,  RDE (\ref{msde}) possesses a unique global solution in the framework of $\beta$-H\"older rough path theory.
\end{proposition}

Next, we will deal with a fast-slow system of RDEs driven by 
mixed fractional Brownian rough path with Hurst index $H\in (\frac13, \frac12]$ of the type
\begin{eqnarray}\label{couple}
\begin{cases}
\begin{aligned} 
X_{t}^{\varepsilon} &=X^\varepsilon_{0}+\int_{0}^{t} f\left(X_{s}^{\varepsilon}, Y_{s}^{\varepsilon}\right) \mathrm{d} s+\int_{0}^{t} \sigma\left(X_{s}^{\varepsilon}\right) \mathrm{d} B_{s}, \\ Y_{t}^{\varepsilon} &=Y^\varepsilon_{0}+\frac{1}{\varepsilon} \int_{0}^{t} g\left(X_{s}^{\varepsilon}, Y_{s}^{\varepsilon}\right) \mathrm{d} s+\frac{1}{\sqrt{\varepsilon}} \int_{0}^{t} h\left(X_{s}^{\varepsilon}, Y_{s}^{\varepsilon}\right) \mathrm{d} W_{s},
\end{aligned}
\end{cases}
\end{eqnarray}
where $X^\varepsilon_{0}=X_0\in \mathbb{R}^m, Y^\varepsilon_{0}=Y_0 \in \mathbb{R}^n$ are arbitrary and non-random but fixed,
while  $\varepsilon$ is a small positive parameter.  The coefficients $f:\mathbb{R}^{m} \times \mathbb{R}^{n} \rightarrow \mathbb{R}^m, g:\mathbb{R}^{m}\times \mathbb{R}^{n} \rightarrow \mathbb{R}^n, \sigma:\mathbb{R}^m \rightarrow \mathbb{R}^{m} \otimes \mathbb{R}^{d}, h:\mathbb{R}^{m} \times \mathbb{R}^{n} \rightarrow \mathbb{R}^{n}\otimes \mathbb{R}^{d'}$ satisfy the following regularity assumptions:
\begin{itemize}
\item (H1) $f$ is a Lipschitz continuous vector field with 
at most linear growth, $\sigma \in C_b^{4}(\mathbb{R}^m; \mathbb{R}^{m}\otimes \mathbb{R}^{d})$.
\item (H2) $g$ is a Lipschitz continuous vector field with at most linear growth, $h \in C_b^{4}(\mathbb{R}^m\times \mathbb{R}^n; \mathbb{R}^{n}\otimes \mathbb{R}^{d'})$.
\end{itemize}
The system of RDEs (\ref{couple}) is a special case of RDE (\ref{msde})
and therefore possesses a unique global solution under (H1) and (H2)
due to Proposition \ref{uniqpro}.

Let us summarize some basic results on the stochastic averaging principle for this kind of slow-fast systems, 
which can be traced back to the work of
Khasminskii \cite{khas1966limit}, see e.g. \cite{duan2014,freidlin2012random,xu2011averaging,xu2014,Xu2015Stochastic,Xu2017Stochastic,Liu2019} and the references therein.  In order to obtain the strong convergence, it is known that the diffusion coefficient $\sigma$ in (\ref{couple}) should not depend on the fast variable $Y^\varepsilon$ (see e.g. \cite{Givon2007}). 
The corresponding literature in the case of perturbation by multiplicative fractional Brownian noise is quite sparse. 
It is worth mentioning that
Hairer and Li \cite{Hairer2019} considered a slow-fast system where the slow component is driven by fBm and proved the convergence to the averaged solution takes place in probability. The most relevant result in our case is the recent work \cite{pei-inahama-xu2020, wangxu2020} which answered affirmatively that an averaging principle still holds for fast-slow mixed SDEs if
disturbances involve both Bm and fBm with $H\in(\frac12,1)$ in the mean square sense. Using the generalized Riemann-Stieltjes integral, an averaging principle in the mean square sense for stochastic partial differential equations driven by fBm subject to an additional fast-varying diffusion process was established in 
\cite{pei-inahama-xu2020b}. We point out that all the aforementioned papers concerning fBm assume $H>\frac12$. 

Therefore, it is quite natural to extend the averaging result to the case $H\in (\frac13, \frac12]$. 
Compared to the known results, 
the main difficulty here is how to deal with the mixed fractional Brownian rough path with Hurst index $H\in (\frac13, \frac12]$. 
We will mainly use Hu and Nualart's fractional calculus approach 
to prove the averaging principle (see Theorem \ref{avethm} below).

Now, following the averaging theory inspired by Khasminskii in \cite{khas1966limit}, we define the averaged RDE as follows:
\begin{eqnarray}\label{xbar}
\bar{X}_t =  \bar X_0+\int_{0}^{t}\bar{f}(\bar X_{s}) \mathrm{d}  s+\int_{0}^{t}\sigma(\bar X_s) \mathrm{d} B_s,
\end{eqnarray}
with $\bar X_0=X_0$, where we set
\begin{eqnarray}\label{fbar}
\bar{f}(\xi)=\int_{\mathbb{R}^{n}} f(\xi, \phi) \mu^{\xi}(\mathrm{d} \phi), \quad \xi \in \mathbb{R}^{m},
\end{eqnarray} 
for a unique invariant measure $\mu^{\xi}$ with respect to the following frozen SDE under condition (H4) below:
\begin{eqnarray}\label{frozon0}
\mathrm{d} Y^{\xi,\phi}_t=\tilde{g}(\xi,Y^{\xi,\phi}_t) \mathrm{d}t+h(\xi,Y^{\xi,\phi}_t)\mathrm{d}^{{\rm I }}{W}_t, \quad Y^{\xi,\phi}_0=\phi  \in \mathbb{R}^{n},
\end{eqnarray}
Here, $\int \cdots \mathrm{d}^{{\rm I }} W$  stands for the usual It\^{o}  integral and
\[
\tilde{g}(\xi,\phi)=g(\xi,\phi)+\frac12\sum_{\bar l=1}^{n} \sum_{\bar j=1}^{d'} \mathcal{D}_h^{(\bar j)} h_{\bar l,\bar j}(\xi,\phi), 
\quad
\mathcal{D}_h^{(\bar j)}=\sum_{\bar l=1}^{n}h_{\bar l, \bar j}(\cdot,\cdot) \partial_{\phi_{\bar l}}.
\]

To establish the averaging principle for  (\ref{couple}), we set the following hypotheses:
\begin{itemize}
	\item  (H3) $f\in C^1_b (\mathbb{R}^m\times \mathbb{R}^n; \mathbb{R}^{m})$.
	\item (H4) There exist $L>0, \beta_{i}>0, i=1,2,$ such that
	\begin{eqnarray*}
		2\langle \phi-\tilde{\phi}, \tilde{g}(\xi, \phi)-\tilde{g}(\xi, \tilde{\phi})\rangle+|h(\xi, \phi)-h(\xi, \tilde{\phi})|^{2} &\leq&-\beta_{1}|\phi-\tilde{\phi}|^{2},\cr
		2\langle \phi, \tilde{g}(\xi, \phi)\rangle+|h(\xi, \phi)|^{2} &\leq&-\beta_{2}|\phi|^{2}+L|\xi|^{2}+L
	\end{eqnarray*}
	for any $\xi \in \mathbb{R}^{m}$ and $\phi, \tilde{\phi} \in \mathbb{R}^{n}$.
\end{itemize}

Now, we present our main result of averaging principle in the $L^1$-sense. 
To our knowledge, this is the first result that proves
the averaging principle for a fast-slow system 
in the framework of rough path theory.
\begin{theorem}\label{avethm}
Let $\frac13 <H\le \frac12$ and $\mathbf{Z}$ be the natural rough path lift of $(B_t, W_t)_{t \in [0,T]}$ as in (\ref{RPZ}).
Assume that $f,\sigma,g,h$ satisfy {\rm (H1)-(H4)}. Then, we have
	\begin{eqnarray*}
\limsup\limits_{\varepsilon \rightarrow 0}\mathbb{E}[\|X^\varepsilon-\bar X\|_{\infty}]  = 0.
	\end{eqnarray*}
Here, $\|\cdot\|_{\infty}$ denotes the supremum norm
over $[0,T]$
and $X^\varepsilon$ and $\bar X$
denote the first level paths of the slow component of
(\ref{couple}) and (\ref{xbar}), respectively. 
\end{theorem}

\begin{remark}
A simple example that satisfies the (H4) is $g(\xi,\phi)=\xi-8\phi$ and $h(\xi,\phi)=\sin \xi +\sin \phi$ when $d =d' =m=n=1$. 
Another example is as follows.
Let $g(\xi,\phi)= - A(\xi) \phi$, where $A$ is a bounded, 
positive, $C^1_b$-function in $\xi$, which is also bounded away from zero.
If  $\|h \|_{\infty} +\|\nabla_{\phi} h \|_{\infty}  +\|\nabla^2_{\phi} h \|_{\infty}$
is sufficiently small, then these $g$ and $h$ satisfy (H4).
Here, $\nabla_{\phi}$ stands for the (partial) gradient with respect to $\phi$
and $\|\cdot\|_{\infty}$ denotes the supremum norm.
 \end{remark}

The rest of the paper is organized as follows. Section \ref{se-2} presents some notations and the pathwise approach based on the techniques of the fractional calculus and rough path theory. The existence and uniqueness theorems to (\ref{msde}) and (\ref{couple}) are proved in Section \ref{se-3}. Section \ref{se-4} is devoted to proving Theorem \ref{avethm}, that is, 
the averaging principle for the fast-slow system driven by mixed fractional Brownian rough path with Hurst index $H\in (\frac13, \frac12]$. 

Throughout this paper, 
 $K$ and $K_\ast$ denote certain positive constants that may vary from line to line. $K_\ast$
is used to emphasize that the constant depends on the corresponding parameter $\ast$, which is one or more than one parameter.

\section{Preliminaries}\label{se-2}
Let $|\cdot|$ stand for an absolute value of a real
number, the Euclidean norm of a finite dimensional vector or of a matrix, $\langle \cdot,\cdot\rangle$ be the Euclidean inner product. Set $\Delta_{a,b}=\{(s,t)~|~a  \leq s\leq t  \leq b\}$. Moreover, for a
function $f:[a,b] \rightarrow \mathbb{R}^m$ we define the following seminorms:
$$\|f\|_{\infty,[a,b]}=\sup_{t\in[a,b]}|f(t)|, \quad \|f\|_{\gamma,[a,b]}=\sup _{(s,t)\in \Delta_{a,b}} \frac{|f(t)-f(s)|}{|t-s|^{\gamma}},$$
where $\gamma\in(0,1]$ and also use the convention $0/0\triangleq0$. We will study continuous $\mathbb{R}^m$-valued paths on some interval $[a,b]$, and we denote the space of such functions by $C([a,b];\mathbb{R}^m)$.
We denote by $C^\gamma([a,b];\mathbb{R}^m)$, the space of $\gamma$-H\"older continuous functions on some interval $[a,b]$ with values in $\mathbb{R}^m$
 and set $\|f\|_{\gamma}:=\|f\|_{\gamma,[0, T]},\|f\|_{\infty}:=\|f\|_{\infty,[0, T]}$  and $\Delta:=\Delta_{0,T}$ for shortness. 
 Next, for a function $v:\Delta_{a,b} \rightarrow \mathbb{R}^m$, that vanishes on
the diagonal, that is, $v(t, t ) = 0$ for $t \in [a, b]$, we set
$$
\|v\|_{\gamma,\Delta_{a,b}}=\sup _{(s,t)\in \Delta_{a,b}} \frac{|v(s,t)|}{|t-s|^{\gamma}}. 
$$
The set of such functions with a finite norm $\|v\|_{\gamma,\Delta_{a,b}}$ is denoted by $ C^{\gamma}_2(\Delta_{a,b};\mathbb{R}^m)$, 
$$
C_{2}^{\gamma}(\Delta_{a,b};\mathbb{R}^m)=\{v:\Delta_{a,b} \rightarrow \mathbb{R}^m~|~v(t, t)=0, t \in[a, b],\|v\|_{\gamma,\Delta_{a,b}}<\infty\}.
$$

The purpose of this paper is to use the pathwise approach including rough path analysis and an approach via fractional calculus to the stochastic calculus with respect to Bm and fBm. To prove Proposition \ref{uniqpro}, we mainly use the rough path analysis developed in Lyons \cite{lyons1994}. To establish the averaging principle (see Theorem \ref{avethm}), we combine the rough path analysis and the pathwise approach via fractional calculus inspired by the work of Hu and Nualart \cite{HuN2009}.

We firstly recall the pathwise approach based on the techniques of the fractional calculus in forthcoming Section 2.1 (See Hu and Nualart \cite{HuN2009} and Definition \ref{defn-frac} below). Then, following the ideas in \cite{Gubinelli2004} and \cite{lyons1994} , we provide an explicit formula for integrals in the rough path sense (see Gubinelli in \cite{Gubinelli2004} and (\ref{riemann}) in Section 2.2).

\subsection{Integrals along rough paths via fractional calculus}
Without approximation by Riemann sums, a pathwise approach to study integration of  vector-valued H\"older continuous rough functions using fractional calculus was firstly developed by Hu and Nualart \cite{HuN2009}. Later,
Ito \cite{YIto2019} showed that the integral defined in Definition \ref{defn-frac} below (see also Definition 3.2 in \cite{HuN2009}) using fractional calculus should be consistent with
that obtained by the usual integration in rough path analysis, that will be given by the limit of the
compensated Riemann-Stieltjes sums in forthcoming (\ref{riemann}) of next subsection. In this subsection, we follow the notations proposed by Hu and Nualart \cite{HuN2009}.

This subsection aims to recall the integrals along rough paths via fractional calculus, that is giving meaning to the following integral 
\begin{eqnarray}\label{oint}
\int^T_0 \sigma(x_r)\mathrm{d}\omega_r =\sum^d_{j=1}\int^T_0 \sigma_j(x_r)\mathrm{d}\omega^j_r, 
\end{eqnarray}
using fractional calculus, where $\omega \in C^\beta([0,T];  \mathbb{R}^d)$ with $\beta\in(\frac13,\frac12)$. Moreover, for a given $\omega \in C^\beta([0,T];  \mathbb{R}^d)$, consider $v$ to be an element of $C_{2}^{2\beta}(\Delta;\mathbb{R}^m\otimes\mathbb{R}^d)$ and assume that the triplets $(x,\omega,v)$ satisfy Chen's relation: for all $0\leq s\le \tau \le t \leq T$ it hold that
\begin{eqnarray}\label{chen2}
v_{s \tau}+v_{\tau t}+(x_{\tau}-x_s) \otimes (\omega_t-\omega_{\tau})=v_{st},
\end{eqnarray}
where $\otimes$ denotes tensor.
\begin{definition}\label{Mmd}
\begin{description}
	\item[(1)] For a given $\omega\in C^\beta([0,T];  \mathbb{R}^d)$,  we denote by $M^\beta_{m,d}$ the space consisting of triplets $(x,\omega,v)\in C^\beta([0,T];  \mathbb{R}^m)\times C^\beta([0,T];  \mathbb{R}^d)\times C_{2}^{2\beta}(\Delta;\mathbb{R}^m\otimes\mathbb{R}^d)$ such that (\ref{chen2}) holds.
	\item[(2)] $(\omega,\omega^2)$ is called a $\beta$-H\"older rough path (over $\mathbb{R}^d$) if $(\omega,\omega,\omega^2)\in M^\beta_{d,d}.$ Moreover, define the space of $\beta$-H\"older rough paths (over $\mathbb{R}^d$) in symbols $\mathscr{C}^\beta([0,T];\mathbb{R}^d)$ and the space of weakly geometric $\beta$-H\"older rough paths in symbols $\mathscr{C}_g^\beta([0,T];\mathbb{R}^d)$ by stipulating that $(\omega,\omega^2)\in \mathscr{C}_g^\beta([0,T];\mathbb{R}^d)$ if and only if $(\omega,\omega^2)\in \mathscr{C}^\beta([0,T];\mathbb{R}^d)$ and $${\rm Sym}(\omega^2_{st})=\frac{1}{2}(\omega_{t}-\omega_{s})\otimes (\omega_{t}-\omega_{s}),$$ for every $s,t.$
\end{description}
\end{definition}

To proceed, for $\alpha\in(0,1)$, define the following fractional derivatives
\begin{eqnarray*}	
	D_{s+}^\alpha f[r]&=&\frac{1}{\Gamma(1-\alpha)}\Big(\frac{f(r)}{(r-s)^\alpha}+\alpha\int_s^r\frac{f(r)-f(\theta)}{(r-\theta)^{\alpha+1}} \mathrm{d} \theta\Big),\cr
	D_{t-}^{1-\alpha} g_{t-} [r]&=&\frac{(-1)^{1-\alpha}}{\Gamma(\alpha)}\Big(\frac{g(r)-g(t)}{(t-r)^{1-\alpha}}+(1-\alpha) \int_r^t\frac{g(r)-g(\theta)}{(\theta-r)^{2-\alpha}}\mathrm{d} \theta\Big),
\end{eqnarray*}	   
where $g_{t-}=g(\cdot)-g(t)$. 

Then, we are now ready to define the integral $\int^b_a \sigma(x_r) \mathrm{d} \omega_r.$
\begin{definition}\label{defn-frac} 
	(cf. \cite[Definition 3.2]{HuN2009}) Let $(x,\omega,v)\in M^\beta_{m,d}$ and $\sigma \in C^{1,\lambda}(\mathbb{R}^m; \mathbb{R}^{m}\otimes \mathbb{R}^{d})$ with $(2+\lambda)\beta>1$. Fix $\alpha \in(0,1)$ such that $1-\beta<\alpha<2\beta$, and $\alpha<\frac{\lambda\beta+1}{2}$. Then, for any $(a,b)\in \Delta$ we define 
	\begin{eqnarray*}
		\int^b_a \sigma(x_r)\mathrm{d}\omega_r&=&(-1)^\alpha \int^b_a\hat{D}^\alpha_{a+} \sigma(x)[r]D^{1-\alpha}_{b-}\omega_{b-}[r]\mathrm{d} r\cr
		&&-(-1)^{2\alpha-1}\int^b_a D^{2\alpha-1}_{a+} \nabla \sigma(x)[r]D^{1-\alpha}_{b-}\mathcal{D}^{1-\alpha}_{b-}v[r]\mathrm{d}r.
	\end{eqnarray*}	
Here, for $r\in (a,b)$,
	\begin{eqnarray*}
		\hat{D}^\alpha_{a+}\sigma(x)[r]=\frac{1}{\Gamma(1-\alpha)}\bigg( \frac{\sigma(x_r)}{(r-a)^\alpha}+\alpha\int^r_a \frac{\sigma(x_r)-\sigma(x_\theta)-\nabla \sigma(x_\theta)(x_r-x_\theta)}{(r-\theta)^{\alpha+1}} \mathrm{d}\theta\bigg)
	\end{eqnarray*}
	is the compensated fractional derivative and 
	\begin{eqnarray*}
		\mathcal{D}_{b-}^{1-\alpha}v[r]=\frac{(-1)^{1-\alpha}}{\Gamma(\alpha)}\bigg(\frac{v_{rb}}{(b-r)^{1-\alpha}}+(1-\alpha) \int_{r}^{b} \frac{v_{rs}}{(s-r)^{2-\alpha}} \mathrm{d} s\bigg)
	\end{eqnarray*}
	is the extension of the fractional derivative of $v$.
\end{definition}
Notice that under the constraints  $1-\beta<\alpha<2\beta$ and $\alpha<\frac{\lambda\beta+1}{2}$, it is easy to prove that the fractional derivatives  $D^{1-\alpha}_{b-}\omega_{b-}[r]$ and $D^{1-\alpha}_{b-}\mathcal{D}^{1-\alpha}_{b-}v[r]$ are well defined because the functions $\omega$ and $\mathcal{D}^{1-\alpha}_{b-}v$ are $\beta$-H\"older continuous. Because there exists a constant $K>0$ such that for all $r,\theta\in[a,b],\theta<r$, we have
\begin{eqnarray*}
|\sigma(x_r)-\sigma(x_\theta)-\nabla \sigma(x_\theta)(x_r-x_\theta)|&\leq& K |r-\theta|^{(1+\lambda)\beta},\cr
|\nabla \sigma(x_r)-\nabla \sigma(x_\theta)| &\leq& K |r-\theta|^{\beta \lambda},
\end{eqnarray*}
then the derivatives $\hat{D}^\alpha_{a+}\sigma(x)[r]$ and $D^{2\alpha-1}_{a+} \nabla \sigma(x)[r]$ are also well-defined. More details can be found e.g. in  \cite[p. 2694]{HuN2009}.

Now, given a continuous path $y\in C([0,T];\mathbb{R}^n)$, we aim to solve the RDE with drift term driven by a $\beta$-H\"older rough path $(\omega,\omega^2)$
\begin{eqnarray}\label{path0}
x_t=x_0+\int_{0}^{t}f(x_s, y_s)\mathrm{d}s+\int_{0}^{t}\sigma(x_s) \mathrm{d}\omega_s,
\end{eqnarray}
where $f\in C^1_b(\mathbb{R}^m\times\mathbb{R}^n; \mathbb{R}^{m}),\sigma\in C_b^3(\mathbb{R}^m; \mathbb{R}^{m}\otimes \mathbb{R}^{d})$. 

The main idea to solve  (\ref{path0}) inspired by Hu and Nualart in \cite[Section 4]{HuN2009} and Garrido-Atienza and Schmalfuss \cite[Theorem 4]{Garrido2018} is to write a system of three components for the enlarged unknown $(x,\omega,v).$ According to Definition \ref{defn-frac}, the first component is just (\ref{path0}) itself, where the right-hand side is a function of $(x, \omega, v)$, i.e.
\begin{eqnarray*}
x_t&=&x_0+\int_{0}^{t}f(x_s, y_s)\mathrm{d}s+(-1)^\alpha \int^t_0\hat{D}^\alpha_{0+} \sigma(x)[r]D^{1-\alpha}_{t-}\omega_{t-}[r]\mathrm{d} r\cr
&&-(-1)^{2\alpha-1}\int^t_0 D^{2\alpha-1}_{0+} \nabla \sigma(x)[r]D^{1-\alpha}_{t-}\mathcal{D}^{1-\alpha}_{t-}v[r]\mathrm{d}r.
\end{eqnarray*}
The second component is
\begin{eqnarray}\label{area0}
v_{st}=\int_{s}^{t}\int_{s}^{r}f(x_q, y_q)\mathrm{d}q \otimes \mathrm{d} \omega_r-\int_{s}^{t}\sigma(x_r) \mathrm{d} \omega^2_{\cdot,t}(r).
\end{eqnarray}	
Note that the second term on the right-hand side of  (\ref{area0}) is a functional of $(x, \omega^2, w)$ again by Definition \ref{defn-frac} ($w$ will be given later), i.e.
\begin{eqnarray}
\int_{s}^{t}\sigma(x_r) \mathrm{d} \omega^2_{\cdot,t}(r)&=&(-1)^\alpha \int^t_s\hat{D}^\alpha_{s+} \sigma(x)[r]D^{1-\alpha}_{t-}\omega^2_{\cdot,t-}[r]\mathrm{d} r\cr
	&&-(-1)^{2\alpha-1}\int^t_s D^{2\alpha-1}_{s+} \nabla \sigma(x)[r]D^{1-\alpha}_{t-}\mathcal{D}^{1-\alpha}_{t-}w_{t,\cdot,t}[r]\mathrm{d}r.
\end{eqnarray}	
The third component is defined by writing $w$ as
a functional of $(x, \omega, v, \omega^2 )$ (see \cite[ (3.26)]{HuN2009}) as follows, for $s \leq q \leq t$
\begin{eqnarray}\label{area00}
w_{t,s,q}
&=&-\frac{(-1)^{2 \alpha-1}}{\Gamma(2-2 \alpha)} \int_{s}^{t}\left[\frac{x_{r}-x_{s}}{(r-s)^{2 \alpha-1}}+(2 \alpha-1) \int_{s}^{r} \frac{x_{r}-x_{\theta}}{(r-\theta)^{2 \alpha}} \mathrm{d} \theta\right]\otimes D_{q-}^{1-\alpha} \mathcal{D}_{q-}^{1-\alpha}\omega^2[r]\mathrm{d} r\cr
&&+\frac{(-1)^{2 \alpha-1}}{\Gamma(2-2 \alpha)} \int_{s}^{q} \left[\frac{\omega_{t}-\omega_{r}}{(r-s)^{2 \alpha-1}}+(2 \alpha-1) \int_{s}^{r} \frac{\omega_{\theta}-\omega_{r}}{(r-\theta)^{2 \alpha}} \mathrm{d} \theta\right] \otimes D_{q-}^{1-\alpha} \mathcal{D}_{q-}^{1-\alpha}v[r] d r \cr
&&-\frac{(-1)^{\alpha}}{\Gamma(1-\alpha)} \int_{s}^{q}\bigg(\frac{(x_{r}-x_{s}) \otimes(\omega_{t}-\omega_{r})}{(r-s)^{\alpha}}+\alpha \int_{s}^{r} \frac{(x_{\theta}-x_{r}) \otimes(\omega_{r}-\omega_{\theta})}{(r-\theta)^{\alpha+1}}\mathrm{d} \theta\bigg)\cr
&& \quad  \otimes D_{q-}^{1-\alpha} \omega_{q-}[r] \mathrm{d}  r. 
\end{eqnarray}

\begin{remark}
It should be noted that the sign in front of 
the second term on the right-hand side of  (\ref{area0}) is negative.
This definition of RDEs was first given in \cite[p. 2701, Eq. (4.2)]{HuN2009} when $f \equiv 0$.
However, a negative sign is missing there (and in many other subsequent works). Concerning this, the right side of (\ref{area00}) and that of \cite[p.2701, Eq.(4.5)]{HuN2009} have the opposite signs. 
Fortunately, 
since what is actually computed is the norm of these terms,
all the results in \cite{HuN2009} remain valid.
\end{remark}

By a slight generalization of Theorem 4 in \cite{Garrido2018}, a solution of  (\ref{path0}) is defined to be an element of $M^\beta_{m,d}$ when (\ref{path0})-(\ref{area00}) have a solution and given a continuous path $y\in C([0,T];\mathbb{R}^n)$ and a rough path $(\omega,\omega^2) \in \mathscr{C}^\beta([0,T]; \mathbb{R}^d)$ with $\beta\in(\frac13,\frac12)$, (\ref{path0}) has a unique global solution under the conditions $f\in C^1_b(\mathbb{R}^m\times\mathbb{R}^n; \mathbb{R}^{m}),\sigma\in C_b^3(\mathbb{R}^m; \mathbb{R}^{m}\otimes \mathbb{R}^{d})$, see Garrido-Atienza and Schmalfuss \cite[Theorem 4]{Garrido2018} for example.

Based on a slight generalization of these results in Hu and Nualart \cite{HuN2009} and Besal\'{u}, Binotto and Rovira \cite{Besalu2020}, it is not difficult to provide quantitative estimates for the integration operator in (\ref{path0}) and (\ref{area0}). 
Let us recall two proposition from \cite{Besalu2020,HuN2009}.  On $M^\beta_{m,d}$, we introduce the following functionals for any $(s,t)\in\Delta$:
\begin{eqnarray}\label{Phi1}
\Phi_{\beta,[s,t]}(x,\omega,v)=\|v\|_{2\beta, \Delta_{s,t}}+\|x\|_{\beta,[s,t]}\|\omega\|_{\beta,[s,t]},
\end{eqnarray}
and
\begin{eqnarray}\label{Phi2}
\Phi_{\beta,[s,t]}(x,\omega,v, \omega^2)=\|\omega\|^2_{\beta,[s,t]}\|x\|_{\beta,[s,t]}+\|\omega\|_{\beta,[s,t]}\|v\|_{2\beta, \Delta_{s,t}}+\|x\|_{\beta,[s,t]}\|\omega^2\|_{2\beta, \Delta_{s,t}}.
\end{eqnarray}

We have the following estimates:
\begin{proposition}\label{fy}
{\rm(cf. \cite[Proposition 3.4]{HuN2009})}
	Let $(x,\omega,v) \in M^\beta_{m,d}$ and $\sigma \in C_b^{1,\lambda}(\mathbb{R}^m; \mathbb{R}^{m}\otimes \mathbb{R}^{d})$ with $(2+\lambda)\beta>1$.
Then, for any $(a,b)\in \Delta$, we have 
	\begin{eqnarray*}
		\bigg|\int_{a}^b \sigma(x_r)\mathrm{d}\omega_r\bigg| &\leq & K|\sigma(x_a)|\|\omega\|_{\beta,[a,b]}(b-a)^\beta\cr
		&&+K \Phi_{\beta,[a,b]}(x,\omega,v)(\|\nabla \sigma\|_\infty +\|\nabla \sigma\|_\lambda\|x\|^\lambda_{\beta,[a,b]}(b-a)^{\lambda\beta})(b-a)^{2\beta},
	\end{eqnarray*}
	where $\Phi_{\beta,[a,b]}(x,\omega,v)$ is defined in (\ref{Phi1}).
\end{proposition}

\begin{proposition}\label{fyz}
{\rm(cf. \cite[Proposition 3.9]{HuN2009})}
	Let $(x,\omega,v) \in M^\beta_{m,d}, (\omega,\omega^2) \in \mathscr{C}^\beta([0,T];\mathbb{R}^d)$ and $\sigma \in C_b^{1,\lambda}(\mathbb{R}^m; \mathbb{R}^{m}\otimes \mathbb{R}^{d})$ with $(2+\lambda)\beta>1$. Then, for any $(a,b)\in \Delta$, we have 
\begin{eqnarray*}
\bigg|\int_{a}^{b}\sigma(x_{r}) \mathrm{d}\omega^2_{\cdot, b}(r)\bigg|& \leq & K|\sigma(x_{a})| \Phi_{\beta, [a, b]}(\omega, \omega, \omega^2)(b-a)^{2 \beta} \cr
&&+K(\|\nabla \sigma\|_{\infty}+\|\nabla \sigma\|_{\lambda}\|x\|_{\beta, [a, b]}^{\lambda}(b-a)^{\lambda \beta}) \Phi_{\beta, [a, b]}(x,\omega,v, \omega^2)(b-a)^{3 \beta},
\end{eqnarray*}
where $\Phi_{\beta,[a,b]}(x,\omega,v, \omega^2)$ is defined in (\ref{Phi2}).
\end{proposition}
Given $(\omega,\omega^2) \in \mathscr{C}^\beta([0,T];\mathbb{R}^d)$, we write
$$\Lambda_\omega=1+\|\omega\|^2_{\beta}+\|\omega^2\|_{2\beta}.$$

In \cite{Besalu2020,Besalu2011}, Besal\'{u} and coauthors derived the upper bound of the supremum norm and H\"older norm of the solution $x$ in (\ref{path0}) without $y$.  By a slight generalization of these results in \cite{Besalu2020,Besalu2011}, it is easy to obtain the following lemma.
\begin{lemma}\label{xestim-om} {\rm (cf. \cite[Theorem 4.1]{Besalu2011} and \cite[Proposition 5.2]{Besalu2020})} 
	Assume that $f\in C^1_b(\mathbb{R}^m\times\mathbb{R}^n; \mathbb{R}^{m}),\sigma\in C_b^3(\mathbb{R}^m; \mathbb{R}^{m}\otimes \mathbb{R}^{d})$. Then, we have the following estimates for the solution $(x,\omega,v) \in M^\beta_{m,d}$ of (\ref{path0}):
	\begin{eqnarray*}
		\|x\|_{\infty}+\|x\|_{\beta} +\|v \|_{2 \beta,\Delta}  \leq K_{\beta, T, |x_0|,f,\sigma}\Lambda_\omega^{\diamond(\beta)},
	\end{eqnarray*}
where $\diamond(\beta)$ denotes a certain positive constant which depends only on $\beta$ and $K_{\beta, T, |x_0|,f,\sigma}>0$ is a constant depending only on $\beta, T, |x_0|,\|f\|_{\infty},\|\sigma\|_{\infty},\|\nabla \sigma\|_{\infty}$ and $\|\nabla^2 \sigma\|_{\infty}$ and is independent of $y$.
\end{lemma}
\para{Proof:}
The proof of this result can be found in \cite{Besalu2020,Besalu2011}, although the paper \cite{Besalu2011} did not deal with any drift term $f$. However, such a term can be handled easily since the boundedness of $f$. Thus, by \cite[Theorem 4.1 (i)]{Besalu2011} and Proposition \ref{fy} and Proposition \ref{fyz}, it is not difficulty to see that 
	\begin{eqnarray*}
	\|x\|_{\infty}  \leq |x_{0}|+T(K \rho_{f, \sigma} \Lambda_{\omega})^{\frac{1}{\beta}}+1, 
\end{eqnarray*}
where $
\rho_{f, \sigma}:=\|f\|_{\infty} T^{1-\beta}+\|\sigma\|_{\infty}+\|\nabla \sigma\|_{\infty}+\|\nabla^2 \sigma\|$ and $K$ is a universal constant depending on $\beta$. 

By \cite[Step 2, p.251-252]{Besalu2011} or \cite[Proposition 5.2]{Besalu2020}, we can find a bound for the H\"older norms of $\|x\|_{\beta}$ and $\|v \|_{2 \beta,\Delta}$, i.e. there exists a constant $K_{\beta, T, |x_0|,f,\sigma}$ depending only on $\beta, T, |x_0|,\|f\|_{\infty},\|\sigma\|_{\infty},\|\nabla \sigma\|_{\infty}$ and $\|\nabla^2 \sigma\|_{\infty}$ such that
	\begin{eqnarray*}
\|x\|_{\beta} +\|v \|_{2 \beta,\Delta}  \leq K_{\beta, T, |x_0|,f,\sigma}\Lambda_\omega^{\diamond(\beta)}.
\end{eqnarray*}
Thus, the statement holds. \qed

In order to give some estimates involving a Lipschitz function $\sigma$, we need to introduce some notation:
\begin{eqnarray*}
	G_{\beta,[a,b]}^{1}(\sigma, x, \tilde{x}, \omega,v)&=& K [(\|\nabla^2 \sigma\|_{\infty}+\|\nabla^2 \sigma\|_{\lambda}(\|x\|_{\beta,[a,b]}^{\lambda}+\|\tilde{x}\|_{\beta,[a,b]}^{\lambda})(b-a)^{\lambda \beta})\cr
	&&\times(\Phi_{\beta,[a,b]}(x, \omega,v)+\|\omega\|_{\beta,[a,b]}\|\tilde{x}\|_{\beta,[a,b]})+\|\omega\|_{\beta,[a,b]}\|\nabla \sigma\|_{\infty}],\cr
	G_{\beta,[a,b]}^{2}(\sigma, x, \tilde{x}, \omega,v)&=& K [\|\nabla^2 \sigma\|_{\infty}(\Phi_{\beta,[a,b]}(x, \omega,v)+\|\omega\|_{\beta,[a,b]}\|\tilde{x}\|_{\beta,[a,b]})(b-a)^{\beta}\cr
	&&+\|\omega\|_{\beta,[a,b]}\|\nabla \sigma\|_{\infty}], \cr
	G_{\beta,[a,b]}^{3}(\sigma, \tilde{x})&=& K [\|\nabla \sigma\|_{\infty}+\|\nabla^2 \sigma\|_{\infty}\|\tilde{x}\|_{\beta,[a,b]}(b-a)^{\beta}].
\end{eqnarray*}
\begin{proposition}\label{fyy}{\rm(cf. \cite[Proposition 6.4]{HuN2009})}
	Let $(x,\omega,v), (\tilde{x},\omega,\tilde{v}) \in M^\beta_{m,d}$. Assume that $\sigma \in C_b^{2,\lambda}(\mathbb{R}^m; \mathbb{R}^{m}\otimes \mathbb{R}^{d})$ with $(2+\lambda)\beta>1$. Then, for any $(a,b)\in \Delta$, we have 
	\begin{eqnarray*}	
		\bigg|\int_{a}^{b}(\sigma(x_{r})-\sigma(\tilde{x}_{r})) \mathrm{d} \omega_{r}\bigg| &\leq& G_{\beta,[a,b]}^{1}(\sigma, x, \tilde{x}, \omega,v)(b-a)^{ \beta}\|x-\tilde{x}\|_{\infty,  [a, b]}\cr
		&&+G_{\beta,[a,b]}^{2}(\sigma, x, \tilde{x}, \omega,v)(b-a)^{2 \beta}\|x-\tilde{x}\|_{\beta,[a,b]}\cr
		&&+G_{\beta,[a,b]}^{3}(\sigma, \tilde{x})(b-a)^{2 \beta}\|v-\tilde{v}\|_{2 \beta,\Delta_{a,b}}.
	\end{eqnarray*}
\end{proposition}

Let us introduce more useful notation:
\begin{eqnarray*}	
	G^4_{\beta,[a,b]}(\sigma,x,\tilde{x},\omega,v,\omega^2)
	&=& K \{ \|\nabla \sigma\|_\infty \Phi_{\beta,[a,b]}(\omega,\omega,\omega^2)\cr
	&&+[ \|\nabla^2 \sigma\|_\infty + \|\nabla^2 \sigma\|_\lambda (\|x\|^\lambda_{\beta,[a,b]}+ \|\tilde{x}\|^\lambda_{\beta,[a,b]}) (b-a)^{\lambda\beta}) \cr
	&& \quad \times ( \Phi_{\beta,[a,b]}(x, \omega,v,\omega^2) + \|\tilde{x}\|_{\beta,[a,b]}
	\Phi_{\beta,[a,b]}(\omega,\omega,\omega^2)]\},\cr
	G^5_{\beta,[a,b]}(\sigma,x,\tilde{x},\omega, v,\omega^2)
	&=& K [ ( \|\nabla \sigma\|_\infty + \|\nabla^2 \sigma\|_\infty \|\tilde{x}\|_{\beta,[a,b]}
	(b-a)^{\beta}) \Phi_{\beta,[a,b]}(\omega,\omega,\omega^2) \cr
	&&+ \|\nabla^2 \sigma\|_\infty \Phi_{\beta,[a,b]}(x,\omega,v,\omega^2) (b-a)^\beta ],\cr
	G^6_{\beta,[a,b]}(\sigma,\tilde{x},\omega)
	&=& K G^3_{\beta,[a,b]}(\sigma,\tilde{x}) \|\omega\|_{\beta,[a,b]}.
\end{eqnarray*}

From the previous results it is possible to prove the following proposition.
\begin{proposition}\label{fyzz}{\rm (cf. \cite[Proposition 4.9]{Besalu2020})}
Let $(x,\omega,v), (\tilde{x},\omega,\tilde{v}) \in M^\beta_{m,d}$ and $(\omega,\omega^2) \in \mathscr{C}^\beta([0,T];\mathbb{R}^d)$. 
	Assume that $\sigma \in C_b^{2,\lambda}(\mathbb{R}^m; \mathbb{R}^{m}\otimes \mathbb{R}^{d})$ with $(2+\lambda)\beta>1$. Then, for any $(a,b)\in \Delta$, we have 
	\begin{eqnarray*}
		\bigg|\int_{a}^{b}(\sigma(x_{r})-\sigma(\tilde{x}_{r}))  \mathrm{d}\omega^2_{\cdot, b}(r)\bigg| &\leq& G_{\beta,[a,b]}^{4}(\sigma,x,\tilde{x},\omega,v,\omega^2)(b-a)^{2\beta}\|x-\tilde{x}\|_{\infty,  [a, b]} \cr
		&&+G_{\beta,[a,b]}^{5}(\sigma,x,\tilde{x},\omega,v,\omega^2)(b-a)^{3 \beta}\|x-\tilde{x}\|_{\beta,[a,b]} \cr
		&&+G_{\beta,[a,b]}^{6}(\sigma, \tilde{x}, \omega)(b-a)^{3 \beta}\|v-\tilde{v} \|_{2 \beta,\Delta_{a,b}}.
	\end{eqnarray*}
\end{proposition}

\subsection{Rough paths theory with approximation by Riemann sums}
In this subsection, we write a $\beta$-H\"older rough path by $(X,X^2)$ instead of $(\omega,\omega^2)$ following the notations by Friz and Hairer \cite{FrizH2017}.
\begin{definition}
We say that a pair $(x,x^{\dagger})$ is a controlled path with respect to $(X,X^2)$
if the following decomposition holds
$$x_t-x_s=x^{\dagger}_s(X_t-X_s)+r_{st}, \quad (s,t) \in \Delta,$$
for certain $x^{\dagger}\in C^{\beta}([0,T];\mathbb{R}^{m}\otimes \mathbb{R}^{d}) $ and $r\in C_2^{2\beta}(\Delta;\mathbb{R}^{m})$ where $x^{\dagger}$ is the Gubinelli derivative of $x$. The totality of such $(x, x^{\dagger})$ is denoted by
$\mathcal{D}_{X}^{2\beta}([0,T];\mathbb{R}^m)$, see e.g. \cite[Definition 4.6]{FrizH2017}. We will omit the value space for simplicity of presentation, i.e. $\mathcal{D}_{X}^{2\beta}$ for shorteness.
\end{definition}
When $(x, x^{\dagger}) \in \mathcal{D}_{X}^{2\beta}$, then
the rough integral of $x$ against $X$ can be defined as
\begin{eqnarray}\label{xXinte}
\int_{s}^{t}x^i_{r}\mathrm{d} X^j_{r}=\lim _{|\mathcal{P}| \rightarrow 0} \sum_{t_{k} \in \mathcal{P}}\big(x^i_{t_{k}}(X^j_{t_{k+1}}-X^j_{t_{k}})+ \sum_{\ell=1}^{d}(x^{\dagger}_{t_{k}})^{i,\ell}X^{2,\ell j}_{t_{k} t_{k+1}}\big),
\end{eqnarray}
where $i=1,\ldots m, j=1,\ldots,d$ and the limit is taken over all partitions $\mathcal{P}=\{s=t_{-1}=t_{0}<t_{1}<\ldots<t_{n}=t_{n+1}=t\}$ such that $|\mathcal{P}|=\sup _{t_{k} \in \mathcal{P}}|t_{k}-t_{k-1}|\rightarrow 0.$ It is known that $(\int_{0}^{\cdot}x^i_r \mathrm{d} X^j_r, x^i  {\mathbf e}_j) \in \mathcal{D}_{X}^{2\beta}$,
where $\{{\mathbf e}_j\}$ is the canonical basis of $\mathbb{R}^d$.

\begin{proposition}\label{varphi}
{\rm (cf. e.g. \cite[Lemma 7.3]{FrizH2017})}  Let $(x, x^{\dagger}) \in \mathcal{D}_{X}^{2\beta}$ and $\varphi \in C_{b}^{2}\left(\mathbb{R}^{m};\mathbb{R}^{n}\right) .$ Then we can define a controlled path $(\varphi(x), D\varphi(x) x^{\dagger})\in \mathcal{D}_{X}^{2\beta}$, that is, $\varphi(x)$ is controlled by $X$ if we take $\varphi(x)^{\dagger}=D\varphi(x) x^{\dagger} $ as a Gubinelli derivative of $\varphi(x)$,  i.e.
\begin{eqnarray*}
\varphi(x_{t})-\varphi(x_{s})-\varphi(x)^{\dagger}_{s}(X_{t}-X_{s})\in C_2^{2\beta}(\Delta;\mathbb{R}^{n}).
\end{eqnarray*}
\end{proposition}
Using appropriate estimates for the integrals, the solution to the following RDE driven by $(X,X^2)$ in the sense of controlled paths theory:
		\begin{eqnarray}\label{12}
	x_t=x_0+\int_{0}^{t}f(x_s)\mathrm{d}s+\int_{0}^{t}\sigma(x_s) \mathrm{d} X_s, 
\end{eqnarray}
with $f: \mathbb{R}^{m} \rightarrow \mathbb{R}^{m}, \sigma: \mathbb{R}^{m} \rightarrow \mathbb{R}^{m}\otimes \mathbb{R}^{d}$ is obtained via a fixed point argument.
\begin{lemma}\label{Gubinelli}
	{\rm (cf. e.g. \cite[Section 8]{FrizH2017})} Suppose that $(X,X^2) \in \mathscr{C}^\beta([0,T];\mathbb{R}^d)$ and let $f\in C_b^{1}(\mathbb{R}^{m}, \mathbb{R}^{m})$ and $\sigma \in C_b^{3}(\mathbb{R}^{m}, \mathbb{R}^{m}\otimes \mathbb{R}^{d})$. Then the equation (\ref{12})
	possesses a unique global solution $(x, \sigma(x)) \in \mathcal{D}_{X}^{2\beta}([0,T],\mathbb{R}^m)$. Here, both sides 
	are understood as elements of $\mathcal{D}_{X}^{2\beta}([0,T],\mathbb{R}^m)$.
\end{lemma}

As a consequence of this Proposition \ref{varphi} and Lemma \ref{Gubinelli}, we have
\begin{eqnarray}\label{riemann}
\int_{0}^{t} \sigma_{j}(x_{r}) \mathrm{d} X_{r}^{j}
=\lim_{|\mathcal{P}| \rightarrow 0} \sum_{t_{k} \in \mathcal{P}}\Big(\sigma_{j}(x_{t_{k}})(X_{t_{k+1}}^{j}-X_{t_{k}}^{j})+\sum_{\ell=1}^{d} \mathcal{D}_{\sigma}^{(\ell)} \sigma_{j}(x_{t_{k}}) X_{t_{k} t_{k+1}}^{2,\ell j}\Big),
\end{eqnarray}
with $\sigma_j: \mathbb{R}^{m} \rightarrow \mathbb{R}^{m}, j=1,2 \ldots d$, where the differential operators $ \mathcal{D}^{(\ell)}_\sigma= \sum_{ l=1}^m \sigma_{l,\ell}(\cdot) \partial_{x_l}$. It is known that  $(x_0+\int_{0}^{\cdot}f(x_s)\mathrm{d}s+\int_{0}^{\cdot}\sigma_{j}(x_r) \mathrm{d} X^j_r, \sigma_{j}(x_\cdot)) \in \mathcal{D}_{X}^{2\beta}([0,T],\mathbb{R}^m).$

\begin{remark} Let $\sigma \in C^{3}_b(\mathbb{R}^m; \mathbb{R}^{m}\otimes \mathbb{R}^{d})$.
	\begin{itemize} 
		\item For $(x, \omega, v) \in M^\beta_{m, d}$, the integral $\int_{0}^{T} \sigma (x_{s}) \mathrm{d}  \omega_{s}$ coincides with the integral defined using the $\frac{1}{\beta}$-variation norm (see Ito \cite[Theorem 2.5 and Remark 2.6]{YIto2019}). This implies that $\int_{0}^{T} \sigma(x_{s}) \mathrm{d}  \omega_{s}$ in Definition \ref{defn-frac} can be given by the limit of Riemann sums of the same form of given in (\ref{riemann}). 
		\item The first level of Hu-Nualart type (Definition \ref{defn-frac}) and Gubinelli type (Definition \ref{Gubinelli}) unique global solution are known to coincide.
	\end{itemize}
\end{remark}

Using the flow method, the solution to the RDE with unbounded drift term is obtaind by Riedel and Scheutzow \cite{Riedel2017}.
\begin{lemma}\label{Riedel2017}
	{\rm (cf. e.g. \cite[Theorem 3.1]{Riedel2017})} Suppose that $(X,X^2) \in \mathscr{C}_g^\beta([0,T];\mathbb{R}^d)$ and assume that $f$ is a locally Lipschitz continuous vector field with at most	linear growth on $\mathbb{R}^{m}$ and $\sigma \in C_b^{4}(\mathbb{R}^m; \mathbb{R}^{m}\otimes \mathbb{R}^{d})$. Then a unique global solution exists for any initial value $x_0$.
\end{lemma}

\section{RDEs driven by mixed fractional Brownian rough path}
\label{se-3}
 We mainly use rough path theory recalled in Section \ref{se-2} to prove that  (\ref{msde}) possesses a unique global solution.  To do that, we recall the following lemma.
\begin{lemma}\label{holder}{\rm(cf. \cite[Lemma 2]{nssaa2015})}
	Suppose $(S(t), \mathcal{F}_t)_{t \in [0,T ]}$ is a stochastic process with $\beta$-H\"older trajectories for
	all $ \beta \in (\frac13, \frac12)$, such that $\mathbb{E}[\|S\|^p_\beta] < \infty$ for all $p \geq 1$. Then, for all $\beta' \in (0, \beta)$, there exist a modification of 
	$(s,t)\rightarrow \int_{s}^{t}(S(u)-S(s))\mathrm{d}^{{\rm I }}W_u$ and an almost surely finite random variable $C_{T,\beta'}$ such that
	$$\bigg|\int_{s}^{t}(S(u)-S(s))\mathrm{d}^{{\rm I }}W_u\bigg|\leq C_{T,\beta'}|t-s|^{\frac12+\beta'},\quad (s,t) \in \Delta.$$
\end{lemma}
Now, we are ready to prove Proposition \ref{uniqpro}.
\para{Proof of Proposition \ref{uniqpro}:}
We understand (\ref{msde}) as following RDE
\begin{eqnarray*}
u_t=u_0+\int_{0}^{t}a(u_s) \mathrm{d}s+\int_{0}^{t}(b,c)(u_s)\mathrm{d}\mathbf{Z}_s,
\end{eqnarray*}
where $\mathbf{Z}=(Z,Z^2)$ is a joint, step-2 rough path lift between the Bm $\mathbf{W}$ and fBm $\mathbf{B}$ (which will be defined below). 
Here $(b,c)$ is the $m\times (d+d')$ block matrix. Set $Z_{st} \triangleq Z_t-Z_s$
and denote
$$
Z_{st}=(B_{st},W_{st})^{{\rm T}}, \quad Z^2_{st}=\bigg(\begin{array}{ll}B^2_{st}& I[B,W]_{st}
\\ I[W, B]_{st}
& W^2_{st} 
\end{array}\bigg),
$$
where $\mathbf{B}=(B, B^2)$ is a canonical geometric rough path (see Friz and Hairer \cite[Section 10.3]{FrizH2017} for example) associated to fBm, $\mathbf{W}=(W, W^2)$ is a geometric rough path in Stratonovich sense (see Friz and Hairer \cite[Section 3]{FrizH2017} for example) associated to Bm, and we define for every $(s,t)\in \Delta$,
\begin{eqnarray*}
	I[B,W]_{st} &\triangleq&  \int_{s}^{t}B_{su} \otimes \mathrm{d}^{{\rm I }} W_u,\cr
	I[W,B]_{st} &\triangleq& W_{st} \otimes B_{st} - \int_{s}^{t} \mathrm{d}^{{\rm I }} W_u  \otimes B_{su}.
\end{eqnarray*}

Firstly, we need to prove the joint, step-2 rough path lift $\mathbf{Z}$ satisfies (\ref{chen2}). For $Z^{2,ij}_{st}$, $i,j\in\{1,\ldots, d\}$ and $i,j\in\{d+1,\ldots, d+d'\}$, it is easy to know that $Z^2$ satisfy the Chen's relation. It remains to demonstrate that when $i\in\{1,\ldots, d\}, j\in\{d+1,\ldots, d+d'\}$ and $j \in\{1,\ldots, d\}, i \in\{d+1,\ldots, d+d'\}$, whether we can obtain same relation. Let us study
$$Z^{2,ij}_{st}=\int_{s}^{t}B^{i}_{su} \mathrm{d}^{{\rm I }} W^j_u, \quad i\in\{1,\ldots, d\}, j \in\{d+1,\ldots, d+d'\},$$
then, if $i \in\{1,\ldots, d\}, j \in\{d+1,\ldots, d+d'\}$, we have for $s\leq u \leq t$,
\begin{eqnarray}\label{chenz}
Z^{2,ij}_{st}-Z^{2,ij}_{su}-Z^{2,ij}_{ut}&=&\int_{s}^{t}B^{i}_{sr}  \mathrm{d} ^{{\rm I }}W^j_r-\int_{s}^{u}B^{i}_{sr} \mathrm{d}^{{\rm I }} W^j_r-\int_{u}^{t}B^{i}_{ur}  \mathrm{d} ^{{\rm I }}W^j_r\cr
&=&\int_{u}^{t}B^{i}_{su} \mathrm{d}^{{\rm I }} W^j_r\cr
&=&{Z}^{i}_{su}{Z}^{j}_{ut}.
\end{eqnarray}

By (\ref{chenz}), it is easy to obtain $Z^{2,ij}_{st}-Z^{2,ij}_{su}-Z^{2,ij}_{ut}
={Z}^{i}_{su}{Z}^{j}_{ut}$ holds for $j \in\{1,\ldots, d\}, i \in\{d+1,\ldots, d+d'\}$.
Thus,  the joint, step-2 rough path lift $\mathbf{Z}$ satisfies (\ref{chen2}). 

 Then, there remains the analytic condition to be checked. By Lemma \ref{holder}, it follows that almost surely
$(Z,Z^2) \in \mathscr{C}^{\beta}([0,T]; \mathbb{R}^{d+d'})$ for any $\beta \in(\frac13, H)$. It is easy to check $Z_{st}^{2,ij}+Z^{2,ji}_{st}=Z_{st}^iZ_{st}^j$. So, $(Z,Z^2) \in \mathscr{C}_g^{\beta}([0,T]; \mathbb{R}^{d+d'})$ for any $\beta \in(\frac13, H)$.

Finally, according to Lemma \ref{Riedel2017}, because $a$ is a locally Lipschitz continuous vector field with at most linear growth on $\mathbb{R}^{e}$ and $b_i,c_j \in C_b^{4}(\mathbb{R}^{e}, \mathbb{R}^{e})~(1\leq i \leq d, 1\leq j \leq d')$. Then,  (\ref{msde})
possesses a unique global solution. By Lemma \ref{Gubinelli}, if $a\in C_b^{1}(\mathbb{R}^{e}, \mathbb{R}^{e})$ and $b_i,c_j \in C_b^{3}(\mathbb{R}^{e}, \mathbb{R}^{e})~(1\leq i \leq d, 1\leq j \leq d')$. 
Then,  (\ref{msde})
possesses a unique global solution.
\qed

Through a similar argument as in the proof of Proposition \ref{uniqpro}, we prove that (\ref{couple}) possesses a unique global solution.
\begin{theorem}
Let $\beta\in(\frac13,H)$ and assume the coefficients of  (\ref{couple}) satisfy {\rm (H1)} and {\rm (H2)}. Then,  (\ref{couple})
	possesses a unique global solution.
\end{theorem}
\para{Proof:} This is just the special case of Theorem 1.1. Let
\begin{eqnarray}\label{uFV}
	u^\varepsilon_t:=\bigg(\begin{array}{c}{X^\varepsilon_t} \\ {Y^\varepsilon_t}\end{array}\bigg),  \quad	F(\xi,\phi):=\bigg(\begin{array}{c}{f(\xi,\phi)} \\ {\frac{1}{\varepsilon}g(\xi,\phi)}\end{array}\bigg) \quad {\rm and}	\quad V(\xi,\phi):=\bigg(
	\begin{matrix}
		\sigma(\xi)  & O \\
		O & \frac{1}{\sqrt{\varepsilon}}h(\xi,\phi)
	\end{matrix}
	\bigg),
\end{eqnarray}
where $F:\mathbb{R}^{m+n}\rightarrow \mathbb{R}^{m+n}$ and $V=(V_1,\ldots,V_{d+d'})$ is $(m+n) \times (d+d')$ matrix-valued.
Then, set $(b,c)=V$ and by (H1) and (H2), (\ref{couple}) possesses a unique global solution. \qed

\begin{remark}
	 It is easy to see that
 $\bar f$ is also a Lipschitz continuous vector field and bounded, because $f$ is a Lipschitz continuous and bounded vector field.
Then, by \cite[Theorem 2.2, Page 121]{Besalu2014}, the averaged RDE (\ref{xbar})
	has a unique global solution. So, the averaged RDE is well-posed in the framework of Hu and Nualart \cite{HuN2009}.
\end{remark}

Now, we study the relation between the fast component of RDE (\ref{couple}) and an It\^o SDE. Note that $\mathbb{E}^B,\mathbb{E}^W$ are the expectation with respect to $B$ and $W$, respectively, so that $\mathbb{E}=\mathbb{E}^B\times \mathbb{E}^W$.
\begin{theorem}\label{fastconv}
	The first level path of the fast component of RDE (\ref{couple}) is the following It\^o SDE
	\begin{eqnarray}\label{euqyito}
	Y^\varepsilon_t=Y_0+\frac{1}{{\varepsilon}}\int_{0}^{t} \tilde{g}(X^\varepsilon_s,Y^\varepsilon_s)  \mathrm{d} s+\frac{1}{\sqrt{\varepsilon}}\int_{0}^{t}h (X^\varepsilon_s,Y^\varepsilon_s) \mathrm{d}^{{\rm I }} {W}_s.
	\end{eqnarray}
	Here, $X^\varepsilon$ is the first level path of the slow component of RDE (\ref{couple}) and $\tilde{g}$ has been defined in (\ref{frozon0}).
\end{theorem}
\para{Proof:}
By (\ref{riemann}), we rewrite the rough path lift terms of the right side of (\ref{couple}) as
\begin{eqnarray}\label{levyarea}
&&\sum_{l=1}^{m+n} \sum_{j=1}^{d+d'}	\int_{0}^{t}{V}_{l,j}(u^\varepsilon_s) \mathrm{d}{Z}^j_s \cr
&&\quad \quad =\lim\limits_{|\mathcal{P}|\rightarrow 0}\sum_{t_k\in \mathcal{P}} \sum_{l=1}^{m+n} \sum_{j=1}^{d+d'}\bigg\{{V}_{l,j}(u^\varepsilon_{t_k}) Z^j_{t_k t_{k+1}}+ \sum_{i=1}^{d+d'} \mathcal{D}_{{V}}^{(i)}{V}_{l,j}(u^\varepsilon_{t_k}) Z^{2,ij}_{t_k t_{k+1}}\bigg\},
\end{eqnarray}
where $\mathcal{D}_{{V}}^{(i)}=\sum_{ l=1}^{m+n}{{V}}_{l,i}(\cdot)\partial_{u_l}, i=1,\ldots,d+d'$.

To prove the theorem, it is sufficient to computate the fast component $Y^\varepsilon$ of (\ref{couple}).
According to (\ref{levyarea}), taking $(m+1) \leq l \leq (m+n)$, we have
\begin{eqnarray*}
&&\lim\limits_{|\mathcal{P}|\rightarrow 0}\sum_{t_k\in \mathcal{P}}\sum_{l=m+1}^{m+n}\sum_{j=1}^{d+d'}\bigg\{{V}_{l,j}(u^\varepsilon_{t_k}) Z^j_{t_k t_{k+1}}+\sum_{i=1}^{d+d'} \mathcal{D}_{{V}}^{(i)}{V}_{l,j}(u^\varepsilon_{t_k}) Z^{2,ij}_{t_k t_{k+1}}\bigg\}\cr
&&\quad \quad = \lim\limits_{|\mathcal{P}|\rightarrow 0}\sum_{t_k\in \mathcal{P}} \bigg\{\sum_{\bar l=1}^{n}
\sum_{\bar j=1}^{d'} h_{\bar l, \bar j}(X^\varepsilon_{t_k},Y^\varepsilon_{t_k}) W^{\bar j}_{t_k t_{k+1}}+\sum_{l=m+1}^{m+n}\sum_{j=d+1}^{d+d'} \sum_{i=1}^{d} \mathcal{D}_V^{(i)} V_{l,j}(u^\varepsilon_{t_k}) Z^{2,ij}_{t_k t_{k+1}}\cr
&&\quad \quad \quad +\sum_{l=m+1}^{m+n} \sum_{j=d+1}^{d+d'} \sum_{i=d+1}^{d+d'} \mathcal{D}_V^{(i)} V_{l,j}(u^\varepsilon_{t_k}) Z^{2,ij}_{t_k t_{k+1}}\bigg\}=:\lim\limits_{|\mathcal{P}|\rightarrow 0} \sum_{i=1}^3\mathbf{I}_i.
\end{eqnarray*}
We will prove that $\lim\limits_{|\mathcal{P}|\rightarrow 0}\mathbf{I}_i$ exists for $i=1,2,3$ in the sense of limit in probability.

For the term $\mathbf{I}_1$, it is easy to obtain
\begin{eqnarray*}
\lim\limits_{|\mathcal{P}|\rightarrow 0}\sum_{t_k\in \mathcal{P}}\sum_{\bar l=1}^{n}
\sum_{\bar j=1}^{d'} h_{\bar l, \bar j}(X^\varepsilon_{t_k},Y^\varepsilon_{t_k}) W^{\bar j}_{t_k t_{k+1}}= \sum_{\bar l=1}^{n}
\sum_{\bar j=1}^{d'} \int_{0}^{t}h_{\bar l, \bar j} (X^\varepsilon_s,Y^\varepsilon_s) \mathrm{d}^{{\rm I }} {W}^{\bar j}_s \quad {\rm in}~L^2,
\end{eqnarray*}
by the definition of the It\^o integral because $h_{\bar l,\bar j}(X^\varepsilon_{s},Y^\varepsilon_{s}), s \in[0,T]$, is continuous, bounded and adapted.

For the second term $\mathbf{I}_2$, set $\mathcal{D}_\sigma^{(\hat i)}=\sum_{\hat l=1}^{m}\sigma_{\hat l, \hat i}(\cdot)\partial_{x_{\hat l}}$, there exists a constant $K>0$, one has
\begin{eqnarray*}
 \mathbb{E}^W[\mathbf{I}^2_2]&=&\mathbb{E}^W\bigg[\bigg(\sum_{t_k\in \mathcal{P}}\sum_{l=m+1}^{m+n}\sum_{j=d+1}^{d+d'} \sum_{i=1}^{d} \mathcal{D}_V^{(i)} V_{l,j}(u^\varepsilon_{t_k}) Z^{2,ij}_{t_k t_{k+1}}\bigg)^2\bigg]\cr
&=&\mathbb{E}^W\bigg[\bigg(\sum_{t_k\in \mathcal{P}}\sum_{\bar l=1}^{n}\sum_{\bar j=1}^{d'} \sum_{\hat i=1}^{d} \mathcal{D}_\sigma^{(\hat i)} h_{\bar l, \bar j}(X^\varepsilon_{t_k},Y^\varepsilon_{t_k}) I[B,W]_{t_k t_{k+1}}^{\hat i,\bar j}\bigg)^2\bigg]\cr
&\leq &K\sum_{\bar l=1}^{n} \sum_{\bar j=1}^{d'}\sum_{\hat  i=1}^{d} \mathbb{E}^W\bigg[\bigg(\sum_{t_k\in \mathcal{P}}\mathcal{D}_\sigma^{(\hat i)}h_{\bar l, \bar j}(X^\varepsilon_{t_k},Y^\varepsilon_{t_k}) \int_{t_k}^{t_{k+1}}(B^{\hat i}_{r}-B^{\hat i}_{t_k}) \mathrm{d} ^{{\rm I }} W^{\bar j}_r \bigg)^2\bigg].
\end{eqnarray*}
To prove this, put $A_{\hat i, \bar l, \bar j}(t):=\mathcal{D}_\sigma^{(\hat i)}h_{\bar l, \bar j}(X^\varepsilon_{t},Y^\varepsilon_{t}), A_{\hat i, \bar l, \bar j}(k):=A_{\hat i, \bar l, \bar j}(t_k)$ and consider 
\begin{eqnarray*}
&&	\mathbb{E}^W\bigg[\bigg(\sum_{k} A_{\hat i, \bar l, \bar j}(k) \int_{t_k}^{t_{k+1}}(B^{\hat i}_{r}-B^{\hat i}_{t_k}) \mathrm{d} ^{{\rm I }} W^{\bar j}_r\bigg)^2  \bigg]\cr
&&\quad \quad =\mathbb{E}^W\bigg[\sum_{k,k'} A_{\hat i, \bar l, \bar j}(k) A_{\hat i, \bar l, \bar j}(k') \bigg(\int_{t_k}^{t_{k+1}}(B^{\hat i}_{r}-B^{\hat i}_{t_k}) \mathrm{d}^{{\rm I }}W^{\bar j}_r\bigg) \bigg(\int_{t_{k'}}^{t_{k'+1}}(B^{\hat i}_{r}-B^{\hat i}_{t_{k'}}) \mathrm{d}^{{\rm I }}W^{\bar j}_r\bigg)  \bigg].
\end{eqnarray*}
If $k<k'$ then $A_{\hat i, \bar l, \bar j}(k) A_{\hat i, \bar l, \bar j}(k') \big(\int_{t_k}^{t_{k+1}}(B^{\hat i}_{r}-B^{\hat i}_{t_k}) \mathrm{d}^{{\rm I }}W^{\bar j}_r\big)$ and $\big(\int_{t_{k'}}^{t_{k'+1}}(B^{\hat i}_{r}-B^{\hat i}_{t_{k'}}) \mathrm{d}^{{\rm I }}W^{\bar j}_r\big)$ are independent so the terms vanish in this case, and similarly if $k>k'$. So we are left with
\begin{eqnarray*}
	\mathbb{E}^W[\mathbf{I}^2_2]	\leq K\sum_{\bar l=1}^{n} \sum_{\bar j=1}^{d'}\sum_{\hat  i=1}^{d}\sum_{t_k\in \mathcal{P}} \mathbb{E}^W\bigg[\bigg(\mathcal{D}_\sigma^{(\hat i)}h_{\bar l, \bar j}(X^\varepsilon_{t_k},Y^\varepsilon_{t_k}) \int_{t_k}^{t_{k+1}}(B^{\hat i}_{r}-B^{\hat i}_{t_k}) \mathrm{d} ^{{\rm I }} W^{\bar j}_r \bigg)^2\bigg].
\end{eqnarray*}
Note that $\mathbb{E}^B\big[\|B^{\hat{i}}\|^2_{\beta}\big]<\infty$, thus, one has  
$\mathbb{E}[\hat{\mathbf{I}}^2_2]= \mathbb{E}^B[\mathbb{E}^W[\mathbf{I}^2_2]] \rightarrow 0$ as $|\mathcal{P}|\rightarrow 0$. 

To proceed, for the third term $\mathbf{I}_3$, one has
\begin{eqnarray*}
\lim\limits_{|\mathcal{P}|\rightarrow 0}	\mathbf{I}_3= \lim\limits_{|\mathcal{P}|\rightarrow 0}\sum_{t_k\in \mathcal{P}}\sum_{\bar l=1}^{n} \sum_{\bar i=1}^{d'}\sum_{\bar j=1}^{d'}\mathcal{D}_h^{(\bar i)} h_{\bar l, \bar j}(X^\varepsilon_{t_k},Y^\varepsilon_{t_k})\big(W^2_{t_k t_{k+1}} \big)^{\bar i,\bar j},
\end{eqnarray*}
where $\mathcal{D}_h^{(\bar i)}=\sum_{\bar l=1}^{n}h_{\bar l, \bar i}(\cdot,\cdot) \partial_{y_{\bar l}}$. Following \cite[Theorem 4]{nssaa2015}, for $\bar i=\bar j$, we have 
\begin{eqnarray*}
	\lim\limits_{|\mathcal{P}|\rightarrow 0}\sum_{t_k\in \mathcal{P}}\sum_{\bar l=1}^{n} \sum_{\bar i=1}^{d'}\mathcal{D}_h^{(\bar i)} h_{\bar l, \bar i}(X^\varepsilon_{t_k},Y^\varepsilon_{t_k})\big(W^2_{t_k t_{k+1}} \big)^{\bar i,\bar i}= \frac12 \sum_{\bar l=1}^{n} \sum_{\bar i=1}^{d'}\int_{0}^{t} \mathcal{D}_h^{(\bar i)} h_{\bar l, \bar i}(X^\varepsilon_{s},Y^\varepsilon_{s}) \mathrm{d} s, \quad {\rm a.s.}
\end{eqnarray*}
Then, by \cite[Theorem 4]{nssaa2015} again, for $\bar i \neq \bar j$, we obtain 
\begin{eqnarray*}
\lim\limits_{|\mathcal{P}|\rightarrow 0} \mathbb{E}\bigg[\bigg(\sum_{t_k\in \mathcal{P}}\sum_{\bar l=1}^{n} \sum_{\bar i=1}^{d'}\sum_{\bar j=1}^{d'}\mathcal{D}_h^{(\bar i)} h_{\bar l, \bar j}(X^\varepsilon_{t_k},Y^\varepsilon_{t_k})\big(W^2_{t_k t_{k+1}} \big)^{\bar i,\bar j}\bigg)^2\bigg] = 0 \quad {\rm in}~L^2.
\end{eqnarray*}
Thus, we have shown (\ref{euqyito}). \qed

\section{Averaging principle}\label{se-4}
In this section, we combine the pathwise approach via fractional calculus and rough path theory to estimate the slow component $X^\varepsilon$ and fast component $Y^\varepsilon$ of RDE (\ref{couple}), respectively. 
Now, let us study the slow component of RDE (\ref{couple}) using the pathwise approach via fractional calculus. By (\ref{uFV}), it is not difficult to show that there exists a triplet $(u^\varepsilon,Z,\tilde{v}^\varepsilon) \in M^\beta_{m+n,d+d'}$ (This section follows the notations proposed by Hu and Nualart \cite{HuN2009}). Note that the slow component $X^\varepsilon$ is the solution of (\ref{path0}) with $(y,\omega)$ replaced by $(Y^\varepsilon, B)$. 
In particular, its \textquotedblleft $v$-component\textquotedblright~is of the following form
\begin{eqnarray}\label{areahat}
v^\varepsilon_{st} = \int_{s}^{t}\int_{s}^{r} f(X^\varepsilon_q, Y^\varepsilon_q) \mathrm{d} q \otimes \mathrm{d} B_r-\int_{s}^{t}\sigma(X^\varepsilon_r) \mathrm{d}B^2_{\cdot,t}(r),
\end{eqnarray}
where the last integral is defined based on fractional calculus theory (see Definition \ref{defn-frac}) and
(\ref{areahat}) is well defined under the conditions (H2) and (H3). The stochastic integral of slow component of (\ref{couple}) is a pathwise integral which depends on $B$ and $B^2$. 

Following the discretization techniques inspired by Khasminskii in \cite{khas1966limit}, we introduce an auxiliary process $(\hat{X}^{\varepsilon},\hat{Y}^{\varepsilon})$ and divide $[0,T]$ into intervals depending of size $\delta<1$, where $\delta$ is a fixed positive number depending on $\varepsilon$ which will be chosen later. Then, we construct $\hat{Y}^{\varepsilon}$ with initial value $\hat{Y}_{0}^{\varepsilon}=Y_0,$
\begin{eqnarray}\label{yhat}
\hat{Y}_{t}^{\varepsilon}=Y_0+\frac{1}{\varepsilon} \int_{0}^{t}\tilde{g}(X_{s(\delta)}^{\varepsilon}, \hat{Y}_{s}^{\varepsilon}) \mathrm{d} s+\frac{1}{\sqrt{\varepsilon}}\int_{0}^{t} h(X_{s (\delta)}^{\varepsilon}, \hat{Y}_{s}^{\varepsilon}) \mathrm{d}^{{\rm I }} W_{s},
\end{eqnarray}
where $s(\delta)=\lfloor {s}{\delta^{-1}} \rfloor \delta$  is the nearest breakpoint preceding $s$. Also, we define the process $\hat{X}^{\varepsilon}$ with initial value $\hat{X}_{0}^{\varepsilon}=X_0,$ by
\begin{eqnarray}\label{xhat}
\hat X^\varepsilon_t=X_0+\int_{0}^{t}f(X_{s(\delta)}^\varepsilon,\hat{Y}^\varepsilon_s) \mathrm{d}  s+\int_{0}^{t}\sigma(X^\varepsilon_s) \mathrm{d} B_s,
\end{eqnarray}
in the rough path sense. Similarly, denote the second component of  (\ref{xhat}) and (\ref{xbar}) as $\hat {v}^\varepsilon$ and $\bar v$, i.e. $(\hat X^\varepsilon,B,\hat v^\varepsilon), (\bar X,B, \bar v) \in M^\beta_{m,d}$, respectively.

\subsection{Some estimates on the solutions $X^\varepsilon$, $\hat X^\varepsilon$, $\bar X$, $Y^\varepsilon$ and $\hat Y^\varepsilon$}
From now on, $\diamond(\beta)$ is a certain positive constant and may vary from line to line. Then, by Lemma \ref{xestim-om}, it is easy to obtain the following lemmas.
\begin{lemma}\label{xestim}
Assume that $f,\sigma$ satisfy {\rm (H1)-(H3)}. Then, we have the following estimates:
	\begin{eqnarray*}
	\|X^\varepsilon\|_{\infty}+\|X^\varepsilon\|_{\beta} +\|v^\varepsilon \|_{2 \beta,\Delta}  \leq K_{\beta, T, |X_0|,f,\sigma}\Lambda_B^{\diamond(\beta)},
\end{eqnarray*}
almost surely, 	where $K_{\beta, T, |X_0|,f,\sigma}>0$ is a constant depending only on $\beta, T, |X_0|,\|f\|_{\infty},\|\sigma\|_{\infty},$ $\|\nabla \sigma\|_{\infty}$ and $\|\nabla^2 \sigma\|_{\infty}$ and $\Lambda_B$ has moments of all order.
\end{lemma}
Using similar techniques, the statements proposed in Lemma \ref{xestim} also hold for $\hat{X}^\varepsilon$ and $\bar X$.

\begin{lemma}\label{x-xh}
Assume that $f,\sigma$ satisfy {\rm (H1)-(H3)}. Then, for $(t,t+\delta)\in \Delta$, we have the following estimate:
\begin{eqnarray*}
\sup_{t\in[0,T]}|X^\varepsilon_{t+\delta}-X^\varepsilon_{t}|  \leq K_{\beta, T, |X_0|,f,\sigma}\Lambda_B^{\diamond(\beta)} \delta^{\beta},
\end{eqnarray*}
almost surely.
\end{lemma}

\begin{lemma}\label{x-xhat}
Assume that $f,\sigma$ satisfy {\rm (H1)-(H3)} and let $(X^\varepsilon,B,v^\varepsilon),(\hat X^\varepsilon,B,\hat{v}^\varepsilon) \in M^{\beta}_{m,d}$ be as in (\ref{couple}) and (\ref{xhat}), respectively. Then, we have the following estimates:
	\begin{eqnarray*}	
		\|X^\varepsilon-\hat X^\varepsilon\|_{\infty}& \leq& K_{\beta, T, |X_0|,f,\sigma}\Lambda_B^{\diamond(\beta)}\delta^{\beta}+	\Phi_\varepsilon,\cr
		\|X^\varepsilon-\hat X^\varepsilon\|_{\beta,[a,b]} &\leq& K_{\beta, T, |X_0|,f,\sigma}\Lambda_B^{\diamond(\beta)} \delta^{\beta}+\Phi_\varepsilon,\cr
		\|v^\varepsilon-\hat v^\varepsilon\|_{2 \beta,\Delta_{a,b}} &\leq& \|B\|_{\beta} ( K_{\beta, T, |X_0|,f,\sigma}\Lambda_B^{\diamond(\beta)} \delta^{\beta}+ \Phi_\varepsilon) ,
	\end{eqnarray*}
almost surely,  where  	
\begin{eqnarray*}	
	\Phi_\varepsilon =\bigg\|\int_{0}^{\cdot}(f(X_{r(\delta)}^\varepsilon,\hat{Y}^\varepsilon_r )-f(X_{r(\delta)}^\varepsilon,Y^\varepsilon_r) )\mathrm{d} r\bigg\|_{\infty}+\bigg\|\int_{0}^{\cdot}(f(X_{r(\delta)}^\varepsilon,\hat{Y}^\varepsilon_r)-f(X_{r(\delta)}^\varepsilon,Y^\varepsilon_r) ) \mathrm{d}r\bigg\|_{\beta}.
	\end{eqnarray*}
\end{lemma}
\para{Proof:} We start studying the supremum norm. By Lemma \ref{x-xh} and (H3), we obtain
\begin{eqnarray*}	
\|X^\varepsilon-\hat X^\varepsilon\|_{\infty}&=&	\bigg\|\int_{0}^{\cdot}(f(X^\varepsilon_s,Y^\varepsilon_s)-f(X_{s(\delta)}^\varepsilon,\hat{Y}^\varepsilon_s) ) \mathrm{d}  s\bigg\|_{\infty}\cr 
&\leq&	\sup_{t\in[0,T]} \int_{0}^{t}|f(X^\varepsilon_s,Y^\varepsilon_s)-f(X_{s(\delta)}^\varepsilon,{Y}^\varepsilon_s) | \mathrm{d}  s\cr 
&&+\bigg\|\int_{0}^{\cdot}(f(X_{s(\delta)}^\varepsilon,\hat{Y}^\varepsilon_s)-f(X_{s(\delta)}^\varepsilon,Y^\varepsilon_s) )\mathrm{d} s\bigg\|_{\infty}\cr
&\leq& K_{\beta, T, |X_0|,f,\sigma}\Lambda_B^{\diamond(\beta)}\delta^{\beta}+	\Phi_\varepsilon, \quad {\rm a.s.}
\end{eqnarray*}

Now, we study the H\"older norm. Using Lemma \ref{x-xh} and by (H3) again, we obtain
\begin{eqnarray*}	
	\|X^\varepsilon-\hat X^\varepsilon\|_{\beta,[a,b]} &=&	\sup_{(s,t)\in \Delta_{a,b} }\frac {|\int_{s}^{t}(f(X^\varepsilon_r,Y^\varepsilon_r)-f(X_{r(\delta)}^\varepsilon,\hat{Y}^\varepsilon_r) ) \mathrm{d}r|}{(t-s)^\beta}\cr
	&\leq& K_{\beta, T, |X_0|,f,\sigma} \Lambda_B^{\diamond(\beta)} (b-a)^{1-\beta}  \delta^{\beta}+\Phi_\varepsilon, \quad {\rm a.s.}
\end{eqnarray*}

Now, we study the H\"older norm $\|v^\varepsilon-\hat v^\varepsilon \|_{2 \beta,\Delta_{a,b}}$. By (\ref{areahat}), Fubini's theorem and the argument proposed in \cite[p.2367]{Garrido2018}, we have
\begin{eqnarray*}	
	\|v^\varepsilon-\hat v^\varepsilon \|_{2 \beta,\Delta_{a,b}}&=& \sup_{(s,t)\in \Delta_{a,b} }\frac{| \int_{s}^{t}\int_{s}^{r} (f(X_{q}^\varepsilon,{Y}^\varepsilon_q) -f(X_{q(\delta)}^\varepsilon,\hat{Y}^\varepsilon_q))\mathrm{d}  q \otimes \mathrm{d} B_r|}{(t-s)^{2\beta}}\cr
	\cr &\leq&  \|B\|_{\beta} (K_{\beta, T, |X_0|,f,\sigma} \Lambda_B^{\diamond(\beta)}(b-a)^{1-\beta}  \delta^{\beta}+ \Phi_\varepsilon) , \quad {\rm a.s.}
\end{eqnarray*}
This completed the proof of Lemma \ref{x-xhat}.
\qed

Next, let us study $\|\hat X^\varepsilon-\bar X\|_{\infty} $.
\begin{lemma}\label{xbar-xhat}
Assume that $f,\sigma$ satisfy {\rm (H1)-(H3)} and let $(\hat X^\varepsilon,B,\hat{v}^\varepsilon), (\bar X,B,\bar{v}) \in M^{\beta}_{m,d}$ be as in (\ref{xhat}) and (\ref{xbar}), respectively. Then, we have the following estimates:
	\begin{eqnarray*}
	\|\hat X^\varepsilon-\bar X\|_{\infty} \leq K_{\beta, T, |X_0|,f,\sigma}2^{\Lambda_B^{\diamond(\beta)}}\Lambda_B^{\diamond(\beta)}(\delta^\beta+\mathbf{A}_1+\mathbf{B}_1+ \Phi_\varepsilon),
	\end{eqnarray*}
\end{lemma}	
almost surely, where
\begin{eqnarray*}
\mathbf{A}_1 &=&\Big\|\int_{0}^{\cdot}(f(X_{s(\delta)}^\varepsilon,\hat{Y}^\varepsilon_s) -\bar{f}(X_{s(\delta)}^\varepsilon)) \mathrm{d}  s\Big\|_{\beta},\cr
\mathbf{B}_1&=& \sup_{(s,t)\in \Delta} \frac{|\int_{s}^{t}\int_{s}^{r}(f(X_{q(\delta)}^\varepsilon,\hat{Y}^\varepsilon_q) -\bar{f}(X_{q(\delta)}^\varepsilon)) \mathrm{d} q \otimes \mathrm{d} B_r|}{(t-s)^{2\beta}}.
\end{eqnarray*}
\para{Proof:} Fix a realization of $(B,B^2)$, then everything in this proof is deterministic. Our first purpose is to estimate the H\"older norm $\|\hat X^\varepsilon-\bar{X}\|_{\beta,[a,b]}$.  By (\ref{xbar}) and (\ref{xhat}), we have  
\begin{eqnarray*}
\|\hat X^\varepsilon-\bar{X}\|_{\beta,[a,b]}
&\leq &\bigg\|\int_{0}^{\cdot}(f(X_{s(\delta)}^\varepsilon,\hat{Y}^\varepsilon_s) -\bar{f}(X_{s(\delta)}^\varepsilon)) \mathrm{d}  s\bigg\|_{\beta}\cr
&&+\bigg\|\int_{a}^{\cdot}(\bar{f}(X_{s(\delta)}^\varepsilon)-\bar{f}(X_{s}^\varepsilon)) \mathrm{d}s\bigg\|_{\beta,[a,b]}\cr
&&+\bigg\|\int_{a}^{\cdot}(\bar{f}(X_{s}^\varepsilon)-\bar{f}(\hat X^\varepsilon_{s})) \mathrm{d}  s\bigg\|_{\beta,[a,b]}\cr
&&+\bigg\|\int_{a}^{\cdot}(\bar{f}(\hat X_{s}^\varepsilon)-\bar{f}(\bar X_{s})) \mathrm{d}  s\bigg\|_{\beta,[a,b]}\cr
&&+\bigg\|\int_{a}^{\cdot}(\sigma(\hat X^\varepsilon_s)-\sigma(\bar X_s)) \mathrm{d} B_s\bigg\|_{\beta,[a,b]}\cr
&&+\bigg\|\int_{a}^{\cdot}(\sigma(X^\varepsilon_s)-\sigma(\hat X^\varepsilon_s)) \mathrm{d} B_s\bigg\|_{\beta,[a,b]}\cr
&=:& \sum_{i=1}^{6} \mathbf{A}_i.
\end{eqnarray*}

Let us study $ \mathbf{A}_2, \mathbf{A}_3$ and $\mathbf{A}_4$. By Lemma \ref{x-xhat}, it is easy to obtain
\begin{eqnarray*}	
\sum_{i=2}^{4} \mathbf{A}_i
&\leq& K (b-a)^{1-\beta}(\sup_{t\in[0,T]}|X_{t}^\varepsilon-X^\varepsilon_{t(\delta)}|+\|X^\varepsilon-\hat X^\varepsilon\|_{\infty}+\|\hat X^\varepsilon-\bar X\|_{\infty,  [a, b]})\cr
&\leq & K(b-a)^{1-\beta}(K_{\beta, T, |X_0|,f,\sigma}\Lambda_B^{\diamond(\beta)}\delta^{\beta}+ \Phi_\varepsilon+\|\hat X^\varepsilon-\bar X\|_{\infty,  [a, b]}).
\end{eqnarray*}

Now we estimate $ \mathbf{A}_5$ and $ \mathbf{A}_6$. Since $\sigma \in C_b^{3}(\mathbb{R}^m; \mathbb{R}^{m}\otimes \mathbb{R}^{d})$, taking $\lambda=1$ in Proposition \ref{fyy},  we have
\begin{eqnarray*}	
\mathbf{A}_5 &\leq& G_{\beta,[a,b]}^{1}(\sigma, \hat X^\varepsilon,\bar X, B,\hat{v}^\varepsilon)\|\hat X^\varepsilon-\bar X\|_{\infty,  [a, b]}\cr
&&+G_{\beta,[a,b]}^{2}(\sigma, \hat X^\varepsilon,\bar X, B,\hat{v}^\varepsilon)(b-a)^{\beta}\|\hat X^\varepsilon-\bar X\|_{\beta,[a,b]}\cr
&&+G_{\beta,[a,b]}^{3}(\sigma, \bar X)(b-a)^{\beta}\|\hat v^\varepsilon-\bar v\|_{2 \beta,\Delta_{a,b}},
\end{eqnarray*}
and by Proposition \ref{fyy}, Lemma \ref{x-xh} and Lemma \ref{x-xhat}, we have
\begin{eqnarray*}	
\mathbf{A}_6 &\leq& G_{\beta,[a,b]}^{1}(\sigma, X^\varepsilon, \hat X^\varepsilon, B,v^\varepsilon)\|X^\varepsilon-\hat X^\varepsilon\|_{\infty,[a,b]}\cr
	&&+G_{\beta,[a,b]}^{2}(\sigma, X^\varepsilon, \hat X^\varepsilon, B,v^\varepsilon)(b-a)^{\beta}\|X^\varepsilon-\hat X^\varepsilon\|_{\beta,[a,b]}\cr
	&&+G_{\beta,[a,b]}^{3}(\sigma, \hat X^\varepsilon)(b-a)^{\beta}\|v^\varepsilon-\hat v^\varepsilon\|_{2 \beta,\Delta_{a,b}}\cr
	&\leq & G_{\beta,[a,b]}^{1}(\sigma, X^\varepsilon, \hat X^\varepsilon, B,v^\varepsilon) ( K_{\beta, T, |X_0|,f,\sigma}\Lambda_B^{\diamond(\beta)} \delta^{\beta} +\Phi_\varepsilon)\cr
	&&+G_{\beta,[a,b]}^{2}(\sigma, X^\varepsilon, \hat X^\varepsilon, B,v^\varepsilon) (b-a)^{\beta}(K_{\beta, T, |X_0|,f,\sigma}\Lambda_B^{\diamond(\beta)} \delta^{\beta}+\Phi_\varepsilon)\cr 
	&&+ G_{\beta,[a,b]}^{3}(\sigma, \hat X^\varepsilon) \|B\|_{\beta} (b-a)^{\beta}( K_{\beta, T, |X_0|,f,\sigma}\Lambda_B^{\diamond(\beta)} \delta^{\beta}+\Phi_\varepsilon).
\end{eqnarray*}

Putting above estimations together, we obtain
\begin{eqnarray*}
\|\hat X^\varepsilon-\bar{X}\|_{\beta,[a,b]}
&\leq & \mathbf{A}_1+  K_{\beta, T, |X_0|,f,\sigma} \Psi_1\Lambda_B^{\diamond(\beta)} \delta^{\beta}+K_{\beta,T} \Psi_1 \Phi_\varepsilon+  \Psi_2 \|\hat X^\varepsilon-\bar X\|_{\infty,  [a, b]}\cr
&&+G_{\beta,[a,b]}^{2}(\sigma, \hat X^\varepsilon, \bar X, B,\hat v^\varepsilon)(b-a)^{\beta}\|\hat X^\varepsilon-\bar X\|_{\beta,[a,b]}\cr
&&+G_{\beta,[a,b]}^{3}(\sigma, \bar X)(b-a)^{\beta}\|\hat v^\varepsilon-\bar v\|_{2 \beta,\Delta_{a,b}},
\end{eqnarray*}
where we set
\begin{eqnarray*}
\Psi_1
	&=&  1
	+ G_{\beta,[0,T]}^{1}(\sigma, X^\varepsilon, \hat X^\varepsilon, B,v^\varepsilon) +G_{\beta,[0,T]}^{2}(\sigma, X^\varepsilon, \hat X^\varepsilon, B, v^\varepsilon)+ G_{\beta,[0,T]}^{3}(\sigma, \hat X^\varepsilon)\|B\|_{\beta},\cr
\Psi_2
&=&K(b-a)^{1-\beta}+ G_{\beta,[0,T]}^{1}(\sigma, \hat X^\varepsilon,\bar X, B,\hat{v}^\varepsilon).
\end{eqnarray*}

Next, by (\ref{areahat}), we have
\begin{eqnarray*}
\|\hat v^\varepsilon-\bar{v}\|_{2\beta,[a,b]}
	&\leq &\sup_{(s,t)\in \Delta} \frac{|\int_{s}^{t}\int_{s}^{q}(f(X_{r(\delta)}^\varepsilon,\hat{Y}^\varepsilon_r) -\bar{f}(X_{r(\delta)}^\varepsilon)) \mathrm{d} r \otimes \mathrm{d}B_{q}|}{(t-s)^{2\beta}}\cr
	&&+\sup_{(s,t)\in \Delta_{a,b}} \frac{|\int_{s}^{t}\int_{s}^{q}(\bar{f}(X_{r(\delta)}^\varepsilon)-\bar{f}(X_{r}^\varepsilon)) \mathrm{d}  r \otimes \mathrm{d}B_{q}|}{(t-s)^{2\beta}}\cr
	&&+\sup_{(s,t)\in \Delta_{a,b}} \frac{|\int_{s}^{t}\int_{s}^{q}(\bar{f}(X_{r}^\varepsilon)-\bar{f}(\hat X^\varepsilon_{r})) \mathrm{d} r \otimes \mathrm{d}B_{q}|}{(t-s)^{2\beta}}\cr
	&&+\sup_{(s,t)\in \Delta_{a,b}} \frac{|\int_{s}^{t}\int_{s}^{q}(\bar{f}(\hat X_{r}^\varepsilon)-\bar{f}(\bar X_{r})) \mathrm{d} r \otimes \mathrm{d}B_{q}|}{(t-s)^{2\beta}}\cr
	&&+\sup_{(s,t)\in \Delta_{a,b}} \frac{|\int_{s}^{t}(\sigma(\hat X^\varepsilon_r)-\sigma(\bar X_r)) \mathrm{d} B^2_{\cdot,t}(r)|}{(t-s)^{2\beta}}\cr
		&&+\sup_{(s,t)\in \Delta_{a,b}} \frac{|\int_{s}^{t}(\sigma(X^\varepsilon_r)-\sigma(\hat X^\varepsilon_r)) \mathrm{d} B^2_{\cdot,t}(r)|}{(t-s)^{2\beta}}\cr
	&=:&\sum_{j=1}^6\mathbf{B}_j.
\end{eqnarray*}

Now, we estimate the norm $\|\hat v^\varepsilon-\bar{v}\|_{2\beta,[a,b]}$. Let us study $ \mathbf{B}_2, \mathbf{B}_3$ and $\mathbf{B}_4$. By Fubini's theorem and the argument proposed in \cite[p.2367]{Garrido2018}, it is easy to obtain
\begin{eqnarray*}	
\sum_{i=2}^{4} \mathbf{B}_i
\leq  \|B\|_{\beta} (K_{\beta, T, |X_0|,f,\sigma}\Lambda_B^{\diamond(\beta)}\delta^{\beta}+K_{\beta,T} \Phi_\varepsilon+K(b-a)^{1-\beta}\|\hat X^\varepsilon-\bar X\|_{\infty,  [a, b]}).
\end{eqnarray*}

Next, thanks to Proposition \ref{fyzz} (taking $\lambda=1$), we have
\begin{eqnarray*}
\mathbf{B}_5 &\leq& G_{\beta,[a,b]}^{4}(\sigma,\hat X^\varepsilon, \bar X,B,\hat v^\varepsilon,B^2)\|\hat X^\varepsilon-\bar X\|_{\infty,  [a, b]} \cr
&&+G_{\beta,[a,b]}^{5}(\sigma,\hat X^\varepsilon, \bar X,B,\hat v^\varepsilon,B^2)(b-a)^{\beta}\|\hat X^\varepsilon-\bar X\|_{\beta,[a,b]} \cr
&&+G_{\beta,[a,b]}^{6}(\sigma, \bar X, B)(b-a)^{\beta}\|\hat v^\varepsilon-\bar v\|_{2 \beta,\Delta_{a,b}},
\end{eqnarray*}
and by Proposition \ref{fyzz} and Lemma \ref{x-xh}, we have
\begin{eqnarray*}
\mathbf{B}_6
	&\leq& G_{\beta,[a,b]}^{4}(\sigma,X^\varepsilon, \hat X^\varepsilon, B,v^\varepsilon,B^2)(K_{\beta, T, |X_0|,f,\sigma}\Lambda_B^{\diamond(\beta)}\delta^{\beta} +\Phi_\varepsilon)\cr
	&&+G_{\beta,[a,b]}^{5}(\sigma,X^\varepsilon, \hat X^\varepsilon, B,v^\varepsilon,B^2)(b-a)^{\beta}(K_{\beta, T, |X_0|,f,\sigma}\Lambda_B^{\diamond(\beta)}\delta^{\beta}+\Phi_\varepsilon)\cr
	&&+G_{\beta,[a,b]}^{6}(\sigma, \hat {X}^\varepsilon, B)(b-a)^{\beta} \|B\|_{\beta}(K_{\beta, T, |X_0|,f,\sigma}\Lambda_B^{\diamond(\beta)}\delta^{\beta}+\Phi_\varepsilon).
\end{eqnarray*}

We take suitable $a$ and $b$ such that
\begin{eqnarray}\label{d1}
	G_{\beta,[0,T]}^{6}(\sigma, \bar X, B)(b-a)^{\beta} \leq \frac{1}{2}.
\end{eqnarray}
So, we can define $\Delta_{\beta}^{1}$ such that
$\Delta_{\beta}^{1}:=(2 G_{\beta,[0,T]}^{6}(\sigma, \bar X, B))^{-\frac{1}{\beta}}$. Then, for $(b-a) \leq \Delta_{\beta}^{1}$ it is easy to obtain
\begin{eqnarray*}
\|\hat v^\varepsilon-\bar v\|_{2 \beta,\Delta_{a,b}} &\leq&
2\mathbf{B}_1+K_{\beta, T, |X_0|,f,\sigma}{\Psi}_{3}\Lambda_B^{\diamond(\beta)}\delta^{\beta}+K_{\beta,T}\Psi_{3} \Phi_\varepsilon+\Psi_{4} \|\hat X^\varepsilon-\bar X\|_{\infty,  [a, b]}
\cr
&&+2G_{\beta,[a,b]}^{5}(\sigma,\hat X^\varepsilon, \bar X,B,\hat v^\varepsilon,B^2)(b-a)^{\beta}\|\hat X^\varepsilon-\bar X\|_{\beta,[a,b]},
\end{eqnarray*}
where
\begin{eqnarray*}
\Psi_{3} &=&2\|B\|_{\beta}+2G_{\beta,[0,T]}^{4}(\sigma, X^\varepsilon, \hat X^\varepsilon, B,v^\varepsilon,B^2)\cr
&&+2G_{\beta,[0,T]}^{5}(\sigma, X^\varepsilon, \hat X^\varepsilon, B,v^\varepsilon,B^2)+2G_{\beta,[0,T]}^{6}(\sigma, \hat X^\varepsilon, B)  \|B\|_{\beta},\cr
\Psi_{4} &=& 2K (b-a)^{1-\beta} \|B\|_{\beta}+2G_{\beta,[0,T]}^{4}(\sigma,\hat X^\varepsilon, \bar X,B,\hat v^\varepsilon,B^2).
\end{eqnarray*}

Next, we have
\begin{eqnarray*}
	\|\hat X^\varepsilon-\bar{X}\|_{\beta,[a,b]}
	&\leq & \mathbf{A}_1+  K_{\beta, T, |X_0|,f,\sigma}\Psi_1 \Lambda_B^{\diamond(\beta)}\delta^{\beta}+K_{\beta,T} \Psi_1 \Phi_\varepsilon+ \Psi_2 \|\hat X^\varepsilon-\bar X\|_{\infty,  [a, b]}\cr
	&&+G_{\beta,[a,b]}^{2}(\sigma,\hat X^\varepsilon, \bar X,B,\hat v^\varepsilon,B^2)(b-a)^{\beta}\|\hat X^\varepsilon-\bar X\|_{\beta,[a,b]}\cr
	&&+G_{\beta,[a,b]}^{3}(\sigma, \bar X)(b-a)^{\beta}\cr
	&& \quad \times [2\mathbf{B}_1+ K_{\beta, T, |X_0|,f,\sigma}\Psi_{3} \Lambda_B^{\diamond(\beta)}\delta^{\beta}+K_{\beta,T} \Psi_3 \Phi_\varepsilon +\Psi_{4} \|\hat X^\varepsilon-\bar X\|_{\infty,  [a, b]}\cr
&&\quad \quad +2G_{\beta,[a,b]}^{5}(\sigma,\hat X^\varepsilon, \bar X,B,\hat v^\varepsilon,B^2)(b-a)^{\beta}\|\hat X^\varepsilon-\bar X\|_{\beta,[a,b]}].
\end{eqnarray*}

Similar to the definition of $\Delta^1_{\beta}$, we take suitable $a$ and $b$ again such that
\begin{eqnarray}\label{d2}
&&G_{\beta,[0,T]}^{2}(\sigma,\hat X^\varepsilon, \bar X,B,\hat v^\varepsilon,B^2)(b-a)^{\beta}\cr
&&\quad \quad +2G_{\beta,[0,T]}^{5}(\sigma,\hat X^\varepsilon, \bar X,B,\hat v^\varepsilon,B^2)(b-a)^{2\beta}G_{\beta,[0,T]}^{3}(\sigma, \bar X)\leq \frac12.
\end{eqnarray}
Then, we have
\begin{eqnarray*}
	\|\hat X^\varepsilon-\bar{X}\|_{\beta,[a,b]}
	&\leq & 2\mathbf{A}_1+ 2 K_{\beta, T, |X_0|,f,\sigma}\Psi_1\Lambda_B^{\diamond(\beta)}\delta^{\beta}+2 K_{\beta,T} \Psi_1 \Phi_\varepsilon+2 \Psi_2 \|\hat X^\varepsilon-\bar X\|_{\infty,  [a, b]}\cr
	&&+2G_{\beta,[a,b]}^{3}(\sigma, \bar X)(b-a)^{\beta}[2\mathbf{B}_1+ K_{\beta, T, |X_0|,f,\sigma}\Psi_{3} \Lambda_B^{\diamond(\beta)}\delta^{\beta}+K_{\beta,T} \Psi_3 \Phi_\varepsilon\cr
	&&+ \Psi_{4} \|\hat X^\varepsilon-\bar X\|_{\infty,  [a, b]}]\cr
	&=& 2\mathbf{A}_1+4G_{\beta,[a,b]}^{3}(\sigma, \bar X)(b-a)^{\beta}\mathbf{B}_1\cr
	&&+ (2\Psi_1+2 G_{\beta,[a,b]}^{3}(\sigma, \bar X)(b-a)^{\beta}\Psi_{3}) (K_{\beta, T, |X_0|,f,\sigma}\Lambda_B^{\diamond(\beta)}\delta^{\beta}+K_{\beta,T} \Phi_\varepsilon) \cr
	&&+(2\Psi_2+2G_{\beta,[a,b]}^{3}(\sigma, \bar X)(b-a)^{\beta}\Psi_{4}) \|\hat X^\varepsilon-\bar X\|_{\infty,  [a, b]}.
\end{eqnarray*}
Putting 
$$
\|\hat X^\varepsilon-\bar X\|_{\infty,[a, b]} \leq |\hat X_a^\varepsilon-\bar X_a|+(b-a)^{\beta}\|\hat X^\varepsilon-\bar X\|_{\beta, [a, b]},
$$
in above equation
we have
\begin{eqnarray*}
\|\hat X^\varepsilon-\bar X\|_{\infty,[a, b]} &\leq& |\hat X_a^\varepsilon-\bar X_a|+(b-a)^{\beta}\big[2\mathbf{A}_1+4G_{\beta,[a,b]}^{3}(\sigma, \bar X)(b-a)^{\beta}\mathbf{B}_1\cr
	&&+ (2\Psi_1+2 G_{\beta,[a,b]}^{3}(\sigma, \bar X)(b-a)^{\beta}\Psi_{3}) (K_{\beta, T, |X_0|,f,\sigma}\Lambda_B^{\diamond(\beta)}\delta^{\beta}+K_{\beta,T} \Phi_\varepsilon) \cr
&&+(2\Psi_2+2G_{\beta,[a,b]}^{3}(\sigma, \bar X)(b-a)^{\beta}\Psi_{4}) \|\hat X^\varepsilon-\bar X\|_{\infty,  [a, b]}\big].
\end{eqnarray*}

Similar to the definition of $\Delta^1_{\beta}$, we take suitable $a$ and $b$ again such that
\begin{eqnarray}\label{d3}
	(b-a)^{\beta}(2\Psi_2+2 G_{\beta,[0,T]}^{3}(\sigma, \bar X)(b-a)^{\beta}\Psi_{4})\leq \frac12,
\end{eqnarray}
and by Lemma \ref{xestim-om} and Lemma \ref{xestim}, it is easy to know 
	\begin{eqnarray*}
		\|\hat X^\varepsilon\|_{\infty}+\|\hat X^\varepsilon\|_{\beta} +\|\hat v^\varepsilon \|_{2 \beta,\Delta}+\|\bar X\|_{\infty}+\|\bar X\|_{\beta} +\|\bar v \|_{2 \beta,\Delta}  \leq K_{\beta, T, |X_0|,f,\sigma}\Lambda_B^{\diamond(\beta)}
	\end{eqnarray*}
holds. Then, we have
\begin{eqnarray*}
	\|\hat X^\varepsilon-\bar X\|_{\infty,[a, b]} \leq
 2|\hat X_a^\varepsilon-\bar X_a|+K_{\beta, T, |X_0|,f,\sigma}\Lambda_B^{\diamond(\beta)}(\delta^\beta+\mathbf{A}_1+\mathbf{B}_1+ \Phi_\varepsilon).
\end{eqnarray*}
Hence, 
\begin{eqnarray}\label{xxx}
\sup_{0\leq t \leq b}|\hat X_t^\varepsilon-\bar X_t| \leq
2\sup_{0\leq t \leq a}|\hat X_t^\varepsilon-\bar X_t|+ K_{\beta, T, |X_0|,f,\sigma}\Lambda_B^{\diamond(\beta)}(\delta^\beta+\mathbf{A}_1+\mathbf{B}_1+ \Phi_\varepsilon).
\end{eqnarray}

Now, we can take suitable $a$ and $b$. There exists $\Delta^{max}_{\beta}$ such that all $a,b$ with $(b-a)\leq\Delta^{max}_{\beta}$ fulfill (\ref{d1}), (\ref{d2}) and (\ref{d3}), then, it is clear that (\ref{xxx}) holds for all $a$ and $b$ such that $(b-a)\leq\Delta^{max}_{\beta}$. Then, choose a certain $M=K_{\beta, T, |X_0|,f,\sigma}\Lambda_B^{\diamond(\beta)}$, we take a partition $0=t_0<t_1<\cdots<t_M=T$ of the interval $[0,T]$ such that $(t_{i+1}-t_i)\leq \Delta^{max}_{\beta}$. Then, 
\begin{eqnarray}
\sup_{0\leq t \leq t_M=T}|\hat X_t^\varepsilon-\bar X_t| &\leq&
2\sup_{0\leq t \leq t_{M-1}}|\hat X_t^\varepsilon-\bar X_t|\cr
&&+
K_{\beta, T, |X_0|,f,\sigma}\Lambda_B^{\diamond(\beta)}(\delta^\beta+\mathbf{A}_1+\mathbf{B}_1+ \Phi_\varepsilon).
\end{eqnarray}
Repeating the process $M$ times we obtain
\begin{eqnarray*}
\sup_{0\leq t \leq T}|\hat X_t^\varepsilon-\bar X_t| &\leq&
2^M |\hat X_0^\varepsilon-\bar X_0|+ \bigg(\sum_{k=0}^{M-1}2^k \bigg) K_{\beta, T, |X_0|,f,\sigma}\Lambda_B^{\diamond(\beta)}(\delta^\beta+\mathbf{A}_1+\mathbf{B}_1+ \Phi_\varepsilon)\cr
&\leq& K_{\beta, T, |X_0|,f,\sigma}2^{\Lambda_B^{\diamond(\beta)}}\Lambda_B^{\diamond(\beta)}(\delta^\beta+\mathbf{A}_1+\mathbf{B}_1+ \Phi_\varepsilon).
\end{eqnarray*}

Finally, we have
\begin{eqnarray*}
	\|\hat X^\varepsilon-\bar X\|_{\infty} \leq K_{\beta, T, |X_0|,f,\sigma}2^{\Lambda_B^{\diamond(\beta)}} \Lambda_B^{\diamond(\beta)}(\delta^\beta+\mathbf{A}_1+\mathbf{B}_1+ \Phi_\varepsilon).
\end{eqnarray*}
This completed the proof. \qed

By Theorem \ref{fastconv}, the first level path of the fast component of RDE (\ref{couple}) is an It\^o SDE (\ref{euqyito}).  Thus, it is easy to derive an upper bound of the supremum norm of the solution $Y^\varepsilon$. 
By \cite[Lemma 4.3, Lemma 4.4]{pei-inahama-xu2020}, the following two lemmas are obtained.
\begin{lemma}\label{ybound}
	Suppose that {\rm (H1)-(H4)} hold.  Then, we have
	$$
	\sup_{t\in [0,T]}\mathbb{E}[|Y^\varepsilon_t|^{2}] \leq K,
	$$
	where $K>0$ is a constant independent of $\varepsilon$.
\end{lemma}
\begin{lemma}\label{yhat} 	
	Suppose that {\rm (H1)-(H4)} hold.  Then, we have
	$$
	\sup_{t\in [0,T]}\mathbb{E}[|Y_{t}^{\varepsilon}-\hat{Y}^\varepsilon_t|^{2}] \leq K \delta,
	$$
	where $K>0$ is a constant independent of $\delta$ and $\varepsilon$.
\end{lemma}
\subsection{Some estimates on the difference between $f$ and $\bar f$}
Now, let us study $ \mathbf{A}_1$. It is easy to see that
\begin{eqnarray*}
	\mathbf{A}_1
	&\leq &\sup_{(s,t)\in \Delta} \bigg\{\frac{ |\int_{s}^{t}(f(X_{r(\delta)}^\varepsilon,\hat{Y}^\varepsilon_r) -\bar{f}(X_{r(\delta)}^\varepsilon)) \mathrm{d}  r|}{(t-s)^{\beta}} \mathbf{1}_{\ell}\bigg\}\cr
	&&+\sup_{(s,t)\in \Delta}\bigg\{\frac{ |\int_{s}^{t}(f(X_{r(\delta)}^\varepsilon,\hat{Y}^\varepsilon_r) -\bar{f}(X_{r(\delta)}^\varepsilon)) \mathrm{d}  r|}{(t-s)^{\beta}} \mathbf{1}_{\ell^c}\bigg\}\cr
	&=:& \mathbf{A}_{11}+ \mathbf{A}_{12},
\end{eqnarray*}	
where $\mathbf{1}_{\cdot}$ is an indicator function, $\ell:=\{t < (\lfloor {s}{\delta^{-1}}\rfloor+2)\delta\}$ and $\ell^c:=\{t \geq (\lfloor {s}{\delta^{-1}}\rfloor+2)\delta\}$.

On the one hand, by (H3) and the fact that $t -s< \lfloor {s}{\delta^{-1}}\rfloor \delta -s + 2\delta \leq  2\delta $,  we have
\begin{eqnarray*}
	\mathbb{E}[\mathbf{A}^2_{11}]
	&\leq&  \mathbb{E}\bigg[\sup_{(s,t)\in \Delta}\bigg\{\frac{|\int_{s}^{t}(f( X_{r(\delta)}^{\varepsilon}, \hat{Y}_{r}^{\varepsilon})-\bar{f}(X_{r(\delta)}^{\varepsilon})) \mathrm{d}r|^2}{(t-s)^{2\beta}}\mathbf{1}_{\ell}\bigg\}\bigg]\cr
	&\leq& K_{\beta,T} \delta.
\end{eqnarray*}
On the other hand, by (H3) and the fact that 
$\lfloor \lambda_1 \rfloor -\lfloor \lambda_2 \rfloor \leq \lambda_1-\lambda_2 +1, $ for $\lambda_1 \geq \lambda_2 \geq 0$, we have
\begin{eqnarray*}
	\mathbb{E}[\mathbf{A}^2_{12}]
	&\leq &K \mathbb{E}\bigg[\sup_{(s,t)\in \Delta}\bigg\{\frac{|\int_{s}^{(\lfloor {s}{\delta^{-1}}\rfloor+1)\delta}(f( X_{r(\delta)}^{\varepsilon}, \hat{Y}_{r}^{\varepsilon})-\bar{f}(X_{r(\delta)}^{\varepsilon})) \mathrm{d}r|^2}{(t-s)^{2\beta}} \mathbf{1}_{\ell^c}\bigg\}\bigg]\cr
		&&+K \mathbb{E}\bigg[\sup_{(s,t)\in \Delta}\bigg\{\frac{|\int^{t}_{ \lfloor {t}{\delta^{-1}}\rfloor \delta}(f( X_{r(\delta)}^{\varepsilon}, \hat{Y}_{r}^{\varepsilon})-\bar{f}(X_{r(\delta)}^{\varepsilon})) \mathrm{d}r|^2}{(t-s)^{2\beta}} \mathbf{1}_{\ell^c}\bigg\}\bigg]\cr
	&&+K\mathbb{E}\bigg[\sup_{(s,t)\in \Delta}\bigg\{\frac{|\int_{(\lfloor {s}{\delta^{-1}}\rfloor+1)\delta}^{\lfloor {t}{\delta^{-1}}\rfloor \delta}(f( X_{r(\delta)}^{\varepsilon}, \hat{Y}_{r}^{\varepsilon})-\bar{f}(X_{r(\delta)}^{\varepsilon})) \mathrm{d}r|^2}{(t-s)^{2\beta}} \mathbf{1}_{\ell^c}\bigg\}\bigg]\cr
	&\leq &K  \mathbb{E}\bigg[\sup_{(s,t)\in \Delta}\bigg\{(t-s)^{1-2\beta}\bigg|\int_{s}^{(\lfloor {s}{\delta^{-1}}\rfloor+1)\delta}(f( X_{r(\delta)}^{\varepsilon}, \hat{Y}_{r}^{\varepsilon})-\bar{f}(X_{r(\delta)}^{\varepsilon})) \mathrm{d}r\bigg| \mathbf{1}_{\ell^c}\bigg\}\bigg]\cr
	&&+K  \mathbb{E}\bigg[\sup_{(s,t)\in \Delta}\bigg\{(t-s)^{1-2\beta}\bigg|\int^{t}_{ \lfloor {t}{\delta^{-1}}\rfloor \delta}(f( X_{r(\delta)}^{\varepsilon}, \hat{Y}_{r}^{\varepsilon})-\bar{f}(X_{r(\delta)}^{\varepsilon})) \mathrm{d}r\bigg| \mathbf{1}_{\ell^c}\bigg\}\bigg]\cr
	&&+K \mathbb{E}\bigg[\sup_{(s,t)\in \Delta}\bigg\{(t-s)^{-2\beta} \bigg|\sum_{k= \lfloor {s}{\delta^{-1}}\rfloor +1}^{\lfloor {t}{\delta^{-1}}\rfloor-1}\int_{k\delta}^{(k+1)\delta}(f( X_{k \delta}^{\varepsilon}, \hat{Y}_{r}^{\varepsilon})-\bar{f}(X_{k \delta}^{\varepsilon})) \mathrm{d}r\bigg|^2\mathbf{1}_{\ell^c}\bigg\}\bigg]\cr
		&\leq &K_{\beta,T} \delta
	+ K\mathbb{E}\bigg[\sup_{(s,t)\in \Delta} \bigg\{(\lfloor {t}{\delta^{-1}}\rfloor-\lfloor {s}{\delta^{-1}}\rfloor-1)(t-s)^{-2\beta}\cr
	&&\quad \times \sum_{ k=\lfloor {s}{\delta^{-1}}\rfloor +1 }^{\lfloor {t}{\delta^{-1}}\rfloor-1} \bigg|\int_{k\delta}^{(k+1)\delta}(f( X_{k \delta}^{\varepsilon}, \hat{Y}_{r}^{\varepsilon})-\bar{f}(X_{k \delta}^{\varepsilon})) \mathrm{d}r\bigg|^2\mathbf{1}_{\ell^c}\bigg\}\bigg]\cr
			&\leq &K_{\beta,T} \delta
	+ K_{\beta,T}\delta^{-1}\mathbb{E}\bigg[ \sum_{ k=0}^{\lfloor {T}{\delta^{-1}}\rfloor-1} \bigg|\int_{k\delta}^{(k+1)\delta}(f( X_{k \delta}^{\varepsilon}, \hat{Y}_{r}^{\varepsilon})-\bar{f}(X_{k \delta}^{\varepsilon})) \mathrm{d}r\bigg|^2\bigg]\cr
	&\leq & K_{\beta,T} \delta+K_{\beta,T} \delta^{-2}\max_{ 0 \leq k \leq\lfloor {T}{\delta^{-1}}\rfloor-1}\mathbb{E}\bigg[\bigg|\int_{k\delta}^{(k+1)\delta}(f(X_{k\delta}^{\varepsilon}, \hat{Y}_{r}^{\varepsilon}) -\bar{f}( X_{k\delta}^{\varepsilon})) \mathrm{d}r\bigg|^2\bigg].
\end{eqnarray*}

Now, by the construction of  $\hat{Y}^\varepsilon$ and a time shift transformation, for any fixed $k$ and $s\in[0,\delta]$, we have
\begin{eqnarray*}
	\hat{Y}_{s+k\delta}^\varepsilon &=&\hat Y_{k\delta}^\varepsilon+\frac{1}{\varepsilon}\int_{k\delta}^{k\delta+s}\tilde{g}(X_{k\delta}^\varepsilon,\hat{Y}_r^\varepsilon)\mathrm{d}r+\frac{1}{\sqrt{\varepsilon}}\int_{k\delta}^{k\delta+s}h(X_{k\delta}^\varepsilon,\hat{Y}_r^\varepsilon)\mathrm{d}^{{\rm I }}W_r\cr
	&=&\hat Y_{k\delta}^\varepsilon+\frac{1}{\varepsilon}\int_0^s\tilde{g}\big(X_{k\delta}^\varepsilon,\hat{Y}_{r+k\delta}^\varepsilon\big) \mathrm{d} r+\frac{1}{\sqrt{\varepsilon}}\int_0^sh\big(X_{k\delta}^\varepsilon,\hat{Y}_{r+k\delta}^\varepsilon\big)\mathrm{d}^{{\rm I }}W_r^{*},
\end{eqnarray*}
where $W_t^{*}=W_{t+k\delta}-W_{k\delta}$ is the shift version of $W_t$, and hence they have the same distribution.

Let $\bar{W}$ be a Bm and independent of $W$. Construct a process $Y^{X_{k\delta}^\varepsilon,\hat Y_{k\delta}^\varepsilon}$ by means of
\begin{eqnarray}\label{yxy}
Y_{s/\varepsilon}^{X_{k\delta}^\varepsilon,\hat Y_{k\delta}^\varepsilon}&=&\hat Y_{k\delta}^\varepsilon+\int_0^{s/\varepsilon}\tilde{g}\big(X_{k\delta}^\varepsilon, Y_{r}^{X_{k\delta}^\varepsilon,\hat {Y}_{k\delta}^\varepsilon}\big)\mathrm{d} r+\int_0^{s/\varepsilon}h\big(X_{k\delta}^\varepsilon, Y_{r}^{X_{k\delta}^\varepsilon,\hat{Y}_{k\delta}^\varepsilon}\big)\mathrm{d}^{{\rm I }}\bar{W}_r\cr
&=&\hat Y_{k\delta}^\varepsilon+\frac{1}{\varepsilon}\int_0^{s}\tilde{g}\big(X_{k\delta}^\varepsilon,Y_{r/\varepsilon}^{X_{k\delta}^\varepsilon,\hat Y_{k\delta}^\varepsilon}\big)\mathrm{d} r+\frac{1}{\sqrt{\varepsilon}}\int_0^{s}h\big(X_{k\delta}^\varepsilon,Y_{r/\varepsilon}^{X_{k\delta}^\varepsilon,\hat Y_{k\delta}^\varepsilon}\big)\mathrm{d}^{{\rm I }}\bar{\bar{W}}^\varepsilon_r,
\end{eqnarray}
where $\bar{\bar{W}}^\varepsilon_t=\sqrt{\varepsilon}\bar{W}_{t/\varepsilon}$ is the scaled version of $\bar{W}_t$. Because both $W^{*}$ and $\bar{\bar{W}}$ are independent of $(X_{k\delta}^\varepsilon,\hat Y_{k\delta}^\varepsilon)$, by comparison, yields
\begin{eqnarray}\label{dist}
(X_{k\delta}^\varepsilon,\{\hat{Y}_{s+k\delta}^\varepsilon\}_{s\in[0,\delta)}\big)\sim\big(X_{k\delta}^\varepsilon,\{Y_{s/\varepsilon}^{X_{k\delta}^\varepsilon,\hat Y_{k\delta}^\varepsilon}\}_{s\in[0,\delta)}),
\end{eqnarray}
where $\sim$ denotes coincidence in distribution sense.
Thus, we have
\begin{eqnarray*}
	\mathbb{E}[\mathbf{A}_{1}^2]
	&\leq & K_{\beta,T} \delta+K_{\beta,T} \delta^{-2}\max_{ 0 \leq k \leq\lfloor {T}{\delta^{-1}}\rfloor-1}\mathbb{E}\bigg[\bigg|\int_{k\delta}^{(k+1)\delta}(f(X_{k\delta}^{\varepsilon}, \hat{Y}_{r}^{\varepsilon}) -\bar{f}( X_{k\delta}^{\varepsilon})) \mathrm{d} r\bigg|^2\bigg]\cr
	&\leq &  K_{\beta,T} \delta+K_{\beta,T} \varepsilon^2 \delta^{-2}\max_{0 \leq k \leq\lfloor {T}{\delta^{-1}}\rfloor-1} \int_{0}^{\frac{\delta}{\varepsilon} }\int_{\theta}^{\frac{\delta}{\varepsilon}}\mathcal{J}_{k}(s,\theta)\mathrm{d} s\mathrm{d}\theta,
\end{eqnarray*}
where
\begin{eqnarray*}
	\mathcal{J}_{k}(s,\theta)
	=\mathbb{E}[\langle f(X^{\varepsilon}_{k\delta},Y^{X^{\varepsilon}_{k\delta},\hat{Y}^{\varepsilon}_{k\delta}}_{s})
	-\bar{f}(X^{\varepsilon}_{k\delta}), 
	f(X^{\varepsilon}_{k\delta},Y^{X^{\varepsilon}_{k\delta},\hat{Y}^{\varepsilon}_{k\delta}}_\theta)
	-\bar{f}(X^{\varepsilon}_{k\delta}) \rangle].
\end{eqnarray*}

Through a similar argument as in \cite[Appendix B]{pei-inahama-xu2020}, i.e., for any $0 \leq \theta \leq s \leq \frac{\delta}{\varepsilon}$ and $k=0,1, \ldots,[T / \delta]-1$, we have
\[
\mathcal{J}_{k}(s,\theta)
\leq K_{T,|X_0|,|Y_0|} e^{-\frac{\beta_1}{2}(s-\theta)}\mathbb{E}[(1+|X^{\varepsilon}_{k\delta}|^{2}+|\hat{Y}^{\varepsilon}_{k\delta}|^{2})]
\leq
K_{T,|X_0|,|Y_0|} e^{-\frac{\beta_1}{2}(s-\theta)},
\]
where $\beta_1$ is defined in (H4). 
Here,  Lemmas \ref{xestim}, \ref{ybound} and \ref{yhat}
were used for the last inequality.

Thus, we have
\begin{eqnarray}\label{A1}
\mathbb{E}[\mathbf{A}_{1}^2]
\leq  K_{\beta,T} (\delta+\varepsilon \delta^{-1}).
\end{eqnarray}

For the term $\mathbf{B}_1$, by Fubini's theorem and the argument proposed in \cite[p.2367]{Garrido2018}, we have
\begin{eqnarray*}
	\mathbf{B}_1
	&=&\sup_{(s,t)\in \Delta}\bigg\{ \frac{|\int_{s}^{t}\int_{s}^{q}(f(X_{r(\delta)}^\varepsilon,\hat{Y}^\varepsilon_r) -\bar{f}(X_{r(\delta)}^\varepsilon)) \mathrm{d} r \otimes \mathrm{d} B_q |}{(t-s)^{2\beta}}\bigg\}\cr
	&\leq&\sup_{(s,t)\in \Delta}\bigg\{ \frac{|\int_{s}^{t}(f(X_{r(\delta)}^\varepsilon,\hat{Y}^\varepsilon_r) -\bar{f}(X_{r(\delta)}^\varepsilon))  \otimes \int_{r}^{t} \mathrm{d} B_q \mathrm{d} r |}{(t-s)^{2\beta}}\bigg\}\cr
	&\leq &\sup_{(s,t)\in \Delta}\bigg\{ \frac{|\int_{s}^{t}(f(X_{r(\delta)}^\varepsilon,\hat{Y}^\varepsilon_r) -\bar{f}(X_{r(\delta)}^\varepsilon))  \otimes \int_{r(\delta)}^{t} \mathrm{d} B_q \mathrm{d} r |}{(t-s)^{2\beta}}\bigg\}\cr
	&&+\sup_{(s,t)\in \Delta}\bigg\{ \frac{|\int_{s}^{t}(f(X_{r(\delta)}^\varepsilon,\hat{Y}^\varepsilon_r) -\bar{f}(X_{r(\delta)}^\varepsilon))  \otimes \int^{r(\delta)}_{r} \mathrm{d} B_q \mathrm{d} r |}{(t-s)^{2\beta}}\bigg\}\cr
	&=:& \mathbf{B}_{11}+ \mathbf{B}_{12}.
\end{eqnarray*}	
Let us study $\mathbf{B}_{11}$. Similarly, we have
\begin{eqnarray*}
	\mathbb{E}[\mathbf{B}^2_{11}]
	&\leq & \mathbb{E}\bigg[\sup_{(s,t)\in \Delta}\bigg\{ \frac{|\int_{s}^{t}(f(X_{r(\delta)}^\varepsilon,\hat{Y}^\varepsilon_r) -\bar{f}(X_{r(\delta)}^\varepsilon))  \otimes \int_{r(\delta)}^{t} \mathrm{d} B_q \mathrm{d} r |^2}{(t-s)^{4\beta}}\mathbf{1}_{\ell}\bigg\}\bigg]\cr
	&&+\mathbb{E}\bigg[\sup_{(s,t)\in \Delta}\bigg\{ \frac{|\int_{s}^{t}(f(X_{r(\delta)}^\varepsilon,\hat{Y}^\varepsilon_r) -\bar{f}(X_{r(\delta)}^\varepsilon))  \otimes \int_{r(\delta)}^{t} \mathrm{d} B_q \mathrm{d} r |^2}{(t-s)^{4\beta}}\mathbf{1}_{\ell^c}\bigg\}\bigg]\cr
	&\leq & K_{\beta,T} \mathbb{E}[\|B\|_{\beta}^2] \delta \cr
	&&+K \mathbb{E}\bigg[\sup_{(s,t)\in \Delta}\bigg\{\frac{|\int_{s}^{(\lfloor {s}{\delta^{-1}}\rfloor+1)\delta}(f( X_{r(\delta)}^{\varepsilon}, \hat{Y}_{r}^{\varepsilon})-\bar{f}(X_{r(\delta)}^{\varepsilon}))  \otimes \int_{r(\delta)}^{t} \mathrm{d} B_q  \mathrm{d}r|^2}{(t-s)^{4 \beta}} \mathbf{1}_{\ell^c}\bigg\}\bigg]\cr
	&&+K \mathbb{E}\bigg[\sup_{(s,t)\in \Delta}\bigg\{\frac{|\int^{t}_{ \lfloor {t}{\delta^{-1}}\rfloor \delta}(f( X_{r(\delta)}^{\varepsilon}, \hat{Y}_{r}^{\varepsilon})-\bar{f}(X_{r(\delta)}^{\varepsilon}))  \otimes \int_{r(\delta)}^{t} \mathrm{d} B_q  \mathrm{d}r|^2}{(t-s)^{4 \beta}} \mathbf{1}_{\ell^c}\bigg\}\bigg]\cr
	&&+K\mathbb{E} \bigg[\sup_{(s,t)\in \Delta}\bigg\{\frac{|\int_{(\lfloor {s}{\delta^{-1}}\rfloor+1)\delta}^{ \lfloor {t}{\delta^{-1}}\rfloor \delta}(f( X_{r(\delta)}^{\varepsilon}, \hat{Y}_{r}^{\varepsilon})-\bar{f}(X_{r(\delta)}^{\varepsilon}))  \otimes \int_{r(\delta)}^{t} \mathrm{d} B_q  \mathrm{d}r|^2}{(t-s)^{4 \beta}} \mathbf{1}_{\ell^c}\bigg\}\bigg]\cr
	&\leq & K_{\beta,T} \mathbb{E}[\|B\|_{\beta}^2] \delta+K  \mathbb{E}\bigg[\sup_{(s,t)\in \Delta}\bigg\{(t-s)^{1-3\beta}\|B\|_{\beta}\cr
	&&\quad \times\bigg|\int_{s}^{(\lfloor {s}{\delta^{-1}}\rfloor+2)\delta}(f( X_{r(\delta)}^{\varepsilon}, \hat{Y}_{r}^{\varepsilon})-\bar{f}(X_{r(\delta)}^{\varepsilon}))  \otimes \int_{r(\delta)}^{t} \mathrm{d} B_q  \mathrm{d}r\bigg| \mathbf{1}_{\ell^c}\bigg\}\bigg]\cr
	&&+K  \mathbb{E}\bigg[\sup_{(s,t)\in \Delta}\bigg\{(t-s)^{1-3\beta}\|B\|_{\beta}\cr
	&&\quad \times\bigg|\int^{t}_{ \lfloor {t}{\delta^{-1}}\rfloor \delta}(f( X_{r(\delta)}^{\varepsilon}, \hat{Y}_{r}^{\varepsilon})-\bar{f}(X_{r(\delta)}^{\varepsilon}))  \otimes \int_{r(\delta)}^{t} \mathrm{d} B_q  \mathrm{d}r\bigg| \mathbf{1}_{\ell^c}\bigg\}\bigg]\cr
	&&+K \mathbb{E}\bigg[\sup_{(s,t)\in \Delta}\bigg\{\frac{\big|\sum_{ \lfloor {s}{\delta^{-1}}\rfloor +1}^{\lfloor {t}{\delta^{-1}}\rfloor-1}\int_{k\delta}^{(k+1)\delta}(f( X_{k \delta}^{\varepsilon}, \hat{Y}_{r}^{\varepsilon})-\bar{f}(X_{k \delta}^{\varepsilon}))  \otimes \int_{k\delta}^{t} \mathrm{d} B_q \mathrm{d}r\big|^2}{(t-s)^{4\beta}}\mathbf{1}_{\ell^c}\bigg\}\bigg]\cr
	&\leq & K_{\beta,T}	\mathbb{E}[\|B\|_{\beta}^2] \delta 
	+ K\mathbb{E}\bigg[\sup_{(s,t)\in \Delta} \bigg\{(\lfloor {t}{\delta^{-1}}\rfloor-\lfloor {s}{\delta^{-1}}\rfloor-1)(t-s)^{-2\beta}	\|B\|_{\beta}^2\cr
	&&\quad \times \sum_{k=\lfloor {s}{\delta^{-1}}\rfloor +1}^{\lfloor {t}{\delta^{-1}}\rfloor-1} \bigg|\int_{k\delta}^{(k+1)\delta}(f( X_{k \delta}^{\varepsilon}, \hat{Y}_{r}^{\varepsilon})-\bar{f}(X_{k \delta}^{\varepsilon})) \mathrm{d}r\bigg|^2\mathbf{1}_{\ell^c}\bigg\}\bigg]\cr
	&\leq &K_{\beta,T}	\mathbb{E}[\|B\|_{\beta}^2] \delta \cr
&&+K_{\beta,T} \delta^{-1} \sum_{ k=0}^{\lfloor {T}{\delta^{-1}}\rfloor-1}\mathbb{E}\bigg[\|B\|_{\beta}^2\bigg|\int_{k\delta}^{(k+1)\delta}(f(X_{k\delta}^{\varepsilon}, \hat{Y}_{r}^{\varepsilon}) -\bar{f}( X_{k\delta}^{\varepsilon})) \mathrm{d}r\bigg|^2\bigg]\cr
	&\leq &K_{\beta,T}	\mathbb{E}[\|B\|_{\beta}^2] \delta +K_{\beta,T} \delta^{-1} \big(\mathbb{E}\big[\|B\|_{\beta}^4\big]\big)^{\frac12}\cr
	&& \times \sum_{ k=0}^{\lfloor {T}{\delta^{-1}}\rfloor-1} \bigg(\mathbb{E}\bigg[\bigg|\int_{k\delta}^{(k+1)\delta}(f(X_{k\delta}^{\varepsilon}, \hat{Y}_{r}^{\varepsilon}) -\bar{f}( X_{k\delta}^{\varepsilon})) \mathrm{d}r\bigg|^4\bigg]\bigg)^{\frac{1}{2}}\cr
		&\leq &K_{\beta,T} \delta +K_{\beta,T} \delta^{-1} \max_{ 0 \leq k \leq\lfloor {T}{\delta^{-1}}\rfloor-1}\bigg(\mathbb{E}\bigg[\bigg|\int_{k\delta}^{(k+1)\delta}(f(X_{k\delta}^{\varepsilon}, \hat{Y}_{r}^{\varepsilon}) -\bar{f}( X_{k\delta}^{\varepsilon})) \mathrm{d}r\bigg|^2 \bigg]\bigg)^{\frac12}\cr
	&\leq & K_{\beta,T} (\delta+\varepsilon^{\frac12} \delta^{-\frac12}).
\end{eqnarray*}	

To proceed, it is easy to obtain 
\begin{eqnarray*}
	\mathbb{E}[\mathbf{B}^2_{12}]
	\leq   K_{\beta,T} \delta^{2\beta}.
\end{eqnarray*}

Thus, we have 
\begin{eqnarray}\label{B1}
	\mathbb{E}[\mathbf{B}_1^2]
	\leq   K_{\beta,T} (\delta^{2\beta}+\varepsilon^{\frac12} \delta^{-\frac12}).
\end{eqnarray}

Now, let us study $\mathbb{E}[\Phi^2_\varepsilon]$. 
\begin{eqnarray*}	
	\mathbb{E}[\Phi^2_\varepsilon] &\leq&2 \mathbb{E} \bigg[\sup_{t\in[0,T]}\bigg|\int_{0}^{t}(f(X_{s(\delta)}^\varepsilon,\hat{Y}^\varepsilon_s )-f(X_{s(\delta)}^\varepsilon,Y^\varepsilon_s) )\mathrm{d} s\bigg|^2\bigg]\cr
	&&+2 \mathbb{E}\bigg[\sup_{(s,t)\in \Delta}\frac {\big|\int_{s}^{t}(f(X_{r(\delta)}^\varepsilon,\hat{Y}^\varepsilon_r)-f(X_{r(\delta)}^\varepsilon,Y^\varepsilon_r)) \mathrm{d}r\big|^2}{(t-s)^{2\beta}}\bigg]\cr
	&=:& 2\mathbf{C}_1+2\mathbf{C}_2.
\end{eqnarray*}

Through a similar argument as in the estimates of	$\mathbb{E}[\mathbf{A}^2_{11}]$ and $\mathbb{E}[\mathbf{A}^2_{12}]$ and by Lemma \ref{yhat}, we have
\begin{eqnarray*}	
	\mathbf{C}_1 &\leq&2 \mathbb{E} \bigg[\sup_{t\in[0,T]}\bigg|\int_{\lfloor {t}{\delta^{-1}}\rfloor}^{t}(f(X_{s(\delta)}^\varepsilon,\hat{Y}^\varepsilon_s )-f(X_{s(\delta)}^\varepsilon,Y^\varepsilon_s) )\mathrm{d} s\bigg|^2\bigg]\cr
	&&+2 \mathbb{E}\bigg[\sup_{t\in[0,T]}\bigg|\sum_{k=0}^{\lfloor {t}{\delta^{-1}}\rfloor-1} \int_{k\delta}^{(k+1)\delta}(f(X_{k \delta}^\varepsilon,\hat{Y}^\varepsilon_s )-f(X_{k\delta}^\varepsilon,Y^\varepsilon_s))\mathrm{d} s\bigg|^2\bigg]\cr
	&\leq& K\delta^2+ K_{\beta,T} \delta^{-2} \max_{ 0 \leq k \leq\lfloor {T}{\delta^{-1}}\rfloor-1}\mathbb{E}\bigg[\bigg| \int_{k\delta}^{(k+1)\delta}(f(X_{k \delta}^\varepsilon,\hat{Y}^\varepsilon_s) -f(X_{k\delta}^\varepsilon,Y^\varepsilon_s))\mathrm{d} s\bigg|^2\bigg]\cr
	&\leq & K\delta^2+  K_{\beta,T} \delta^{-1}\max_{ 0 \leq k \leq\lfloor {T}{\delta^{-1}}\rfloor-1}\int_{k\delta}^{(k+1)\delta}\mathbb{E}[|f(X_{k \delta}^\varepsilon,\hat{Y}^\varepsilon_s )-f(X_{k \delta}^\varepsilon,Y^\varepsilon_s)|^2] \mathrm{d} s \cr
	&\leq & K_{\beta,T} \delta,
\end{eqnarray*}
and
\begin{eqnarray*}	
	\mathbf{C}_2
	&\leq & K_{\beta,T} \delta+K_{\beta,T} \delta^{-2}\max_{ 0 \leq k \leq\lfloor {T}{\delta^{-1}}\rfloor-1}\mathbb{E}\bigg[\bigg|\int_{k\delta}^{(k+1)\delta}(f(X_{k \delta}^\varepsilon,\hat{Y}^\varepsilon_s)-f(X_{k \delta}^\varepsilon,Y^\varepsilon_s)) \mathrm{d} s\bigg|^2 \bigg]\cr
	&\leq & K_{\beta,T} \delta+K_{\beta,T} \delta^{-1}\max_{ 0 \leq k \leq\lfloor {T}{\delta^{-1}}\rfloor-1}\int_{k\delta}^{(k+1)\delta}\mathbb{E}[|f(X_{k \delta}^\varepsilon,\hat{Y}^\varepsilon_s) )-f(X_{k \delta}^\varepsilon,Y^\varepsilon_s) )|^2] \mathrm{d} s \cr
	&\leq & K_{\beta,T} \delta.
\end{eqnarray*}

Thus, we have
\begin{eqnarray}\label{Psi}
\mathbb{E}[\Phi^2_\varepsilon]
\leq  K_{\beta,T} \delta.
\end{eqnarray}
\subsection{Proof of Theorem \ref{avethm}}
By Lemma \ref{x-xhat} and Lemma \ref{xbar-xhat}, it is easy to have
\begin{eqnarray*}
\|X^\varepsilon-\bar X\|_{\infty}&\leq&\|X^\varepsilon-\hat X^\varepsilon\|_{\infty}+\|\hat X^\varepsilon-\bar X\|_{\infty}\cr
	&\leq& K_{\beta, T, |X_0|,f,\sigma}2^{\Lambda_B^{\diamond(\beta)}}\Lambda_B^{\diamond(\beta)}(\delta^\beta+ \mathbf{A}_1+\mathbf{B}_1+ \Phi_\varepsilon), \quad{\rm a.s.}
\end{eqnarray*}
Next, for each $R>1$, set $D :=\{\Lambda_B \leq R\}$ and $D^c :=\{\Lambda_B > R\}$, 
then, by (\ref{A1}), (\ref{B1}) and (\ref{Psi}), we have
\begin{eqnarray}\label{xx1}
	\mathbb{E}[\|X^\varepsilon-\bar X\|_{\infty}\mathbf{1}_D] 
	&\leq&\mathbb{E}[K_{\beta, T, |X_0|,f,\sigma}2^{\Lambda_B^{\diamond(\beta)}} \Lambda_B^{\diamond(\beta)}(\delta^\beta+ \mathbf{A}_1+\mathbf{B}_1+ \Phi_\varepsilon)] \cr
	&\leq &K_R K_{\beta, T, |X_0|,f,\sigma}(\delta^\beta +\varepsilon^{\frac14} \delta^{-\frac14}),
\end{eqnarray}
where $K_R>0$ is a constant and
\begin{eqnarray}\label{xx2}
	\mathbb{E}[\|X^\varepsilon-\bar X\|_{\infty}\mathbf{1}_{D^c}] 
	&\leq& (\mathbb{E}[\|X^\varepsilon-\bar X\|^2_{\infty}])^{\frac12} P(\Lambda_B > R)^{ \frac12}\cr
		&\leq& (\mathbb{E}[\|X^\varepsilon\|^2_{\infty}+\|\bar X\|^2_{\infty}])^{\frac12} P(\Lambda_B > R)^{ \frac12} \cr
		&\leq& K_{\beta, T, |X_0|,f,\sigma}(\mathbb{E}[\Lambda_B^{\diamond(\beta)}])^{\frac12}P(\Lambda_B > R)^{\frac12},
\end{eqnarray}
where $K_{\beta, T, |X_0|,f,\sigma}$, $\mathbb{E}[\Lambda_B^{\diamond(\beta)}]$ and $P(\Lambda_B > R)$ are all independent of $\varepsilon$.

Putting (\ref{xx1}) and (\ref{xx2}) together and choose $\delta:=\varepsilon\sqrt{-\ln \varepsilon}$, we have 
\begin{eqnarray*}
\limsup\limits_{\varepsilon \rightarrow 0}
 \mathbb{E}[\|X^\varepsilon-\bar X\|_{\infty}] 
 \leq K P(\Lambda_B > R)^{\frac12},
\end{eqnarray*}
where $K>0$ is a constant which is independent of $\varepsilon$ and $R$. Then, let $R \rightarrow \infty$, we have 
\begin{eqnarray*}
\limsup\limits_{\varepsilon \rightarrow 0}
\mathbb{E}[\|X^\varepsilon-\bar X\|_{\infty}] =0.
\end{eqnarray*}
Thus,  the statement of Theorem \ref{avethm} is obtained.

\section*{Acknowledgement}
\quad The authors are grateful to Professor Yu Ito for stimulating discussions during the preparation
of this work. They also thank Professors David Nualart and Yaozhong Hu for helpful comments.
B. Pei would like to thank JSPS for Postdoctoral Fellowships for Research in Japan (Standard).

B. Pei was partially supported by the National Natural Science Foundation of China under
Grant No. 11802216, the Fundamental Research Funds for the Central Universities, the Young
Talent fund of University Association for Science and Technology in Shaanxi, China, and JSPS
Grant-in-Aid for JSPS Fellows under Grant No. JP18F18314. Y. Inahama was partially supported by JSPS KAKENHI under Grant No. JP20H01807 and Grant-in-Aid for JSPS Fellows
under Grant No. JP18F18314. Y. Xu was partially supported by the National Natural Science
Foundation of China under Grant No. 12072264, the Research Funds for Interdisciplinary Subject of Northwestern Polytechnical University, the Shaanxi Provincial Key R\&D Program under
Grant No. 2019TD-010 and No. 2020KW-013.

\bigskip
\begin{flushleft}
	\begin{tabular}{ll}
		Bin \textsc{Pei}
		\\
School of Mathematical Sciences, Fudan University,
\\
220 Handan Rd., Yangpu District, Shanghai, 200433, China
		\\
		\quad and 
		\\
		Faculty of Mathematics,
		Kyushu University,
		\\
		744 Motooka, Nishi-ku, Fukuoka, 819-0395, Japan.    
		\\
		Email: {\tt binpei@hotmail.com}
	\end{tabular}
\end{flushleft}

\begin{flushleft}
	\begin{tabular}{ll}
		Yuzuru \textsc{Inahama}
		\\
		Faculty of Mathematics,
		Kyushu University,
		\\
		744 Motooka, Nishi-ku, Fukuoka, 819-0395, Japan.
		\\
		Email: {\tt inahama@math.kyushu-u.ac.jp}
	\end{tabular}
\end{flushleft}

\begin{flushleft}
	\begin{tabular}{ll}
		Yong \textsc{Xu}
		\\
School of Mathematics and Statistics, Northwestern Polytechnical University,
		\\
	127 West Youyi Rd., Beilin District, Xi'an, Shaanxi, 710072, China.
		\\
		Email: {\tt hsux3@nwpu.edu.cn}
	\end{tabular}
\end{flushleft}
\end{document}